\documentclass{article}

\usepackage{amsmath}
\usepackage{fullpage}
\usepackage{amssymb}
\usepackage{amsfonts}
\usepackage{tikz}
\usepackage{subcaption}
\usepackage{float}
\usepackage{authblk}

\usepackage[square,comma,numbers,sort&compress]{natbib}
\usepackage{hyperref}

\usetikzlibrary{arrows.meta}
\usetikzlibrary{calc}
\usetikzlibrary{shapes.misc}

\renewcommand{\Im}{\operatorname{Im}}

\DeclareMathOperator{\sech}{sech}
\renewcommand{\Re}{\operatorname{Re}}

\renewcommand{\i}{\mathrm{i}}
\newcommand{\e}{\mathrm{e}}
\newcommand{\eps}{\varepsilon}
\newcommand{\as}{\mathrm{as}}

\newcommand{\pdiff}[2]{\frac{\partial #1}{\partial #2}}

\renewcommand{\d}{ \mathrm{d}}
\newcommand{\diff}[2]{\frac{\mathrm{d} #1}{\mathrm{d} #2}}

\title{The emergence of rapid oscillations in a fifth-order KdV equation}
\author[1]{Christopher J. Lustri\footnote{Corresponding author: \texttt{christopher.lustri@sydney.edu.au}}}
\author[2]{S. Jonathan Chapman}
\author[3]{Richard C. P. Nicotra}
\affil[1,3]{School of Mathematics and Statistics, The University of Sydney, NSW 2006, Australia}
\affil[2]{Mathematical Institute, University of Oxford, Andrew Wiles Building, Radcliffe Observatory Quarter, Woodstock Road, Oxford OX2 6GG, UK
}

\begin{document}
\maketitle

\begin{abstract}

We consider the emergence of rapidly propagating, exponentially small-amplitude oscillations in the impulsively started fifth-order Korteweg--de Vries equation (KdV5), and demonstrate that these waves originate due to Stokes' phenomenon in a small-time boundary layer. We rescale the problem in time and construct an auxiliary linear PDE (lKdV5) that reproduces the relevant exponential asymptotic structure, and proceed by two main methods. Firstly, the Fourier transform solution of the auxiliary PDE is analysed via the method of steepest descents to extract the relevant asymptotic contributions and the locations of the Stokes curves. Secondly, we reproduce these results using a direct factorial-over-power analysis. In this analysis, we make use of the initial and far-field conditions to avoid performing the local asymptotic matching that is typically required by exponential asymptotic methods and is frequently impossible for nonlinear PDEs. A similar factorial-over-power analysis is then applied to the full nonlinear PDE, which cannot be studied using steepest descent methods. The results agree well with numerical simulations, demonstrating both the rapid propagation of non-decaying oscillations and transient ripples whose amplitudes decay in the far field, and revealing the complex Stokes structure (including both normal and higher-order Stokes curves) that governs their emergence. Such a complex Stokes structure in a temporal boundary layer is a novel phenomenon.

\end{abstract}

\section{Introduction and Background}

The fifth-order Korteweg--de Vries (KdV) equation, sometimes known as the Kawahara equation, has been used as a model for small-amplitude water waves with Bond number less than, but close to, $1/3$ and Froude number near 1 \cite{hunter1988existence}. It is given by 
\begin{equation}\label{e:kdv50}
\eps^2 w_{xxxxx} + w_{xxx} + 6 w w_x + w_t = 0,
\end{equation}
where the regime of interest is $0 < \epsilon \ll 1$. This equation is a singular perturbation of the integrable KdV equation,
\begin{equation}\label{e:kdv}
w_{xxx} + 6 w w_x + w_t = 0,
\end{equation}
which possesses the solitary-wave solution
\begin{equation}\label{e:kdvsoliton}
w(x,t) = \tfrac{c}{2}\sech^2\left(\tfrac{\sqrt{c}}{2}\xi\right), \qquad \xi = x - c t.
\end{equation}

 The fifth-order KdV equation \ref{e:kdv50} is typically derived for free-surface flows by beginning with the linear shallow water wave equations and retaining higher-order corrections to the dispersive effects. Hunter and Scheurle \cite{hunter1988existence} showed that any solitary-wave solution of \eqref{e:kdv50} must possess non-decaying oscillatory tails, while Pomeau \emph{et al.} \cite{pomeau1988structural} demonstrated that these oscillations have exponentially small amplitude in an appropriate asymptotic limit. As a result, the existence and behaviour of solitary waves of \eqref{e:kdv50} became a canonical problem in exponential asymptotics, with many subsequent studies on related singularly perturbed KdV equations \cite{akylas1992solitary,amick1991singular,benilov1993generation,boyd1991weakly,calvo2000stability,sun1991existence} investigating properties of these nonlocal (or generalised) solitary waves. 

While most of these studies have concentrated on the stationary version of the equation, there has been some work on the non-stationary case. The authors of \cite{benilov1993generation} performed numerical simulations of solutions to \eqref{e:kdv50} with an initial pulse. They found that the oscillations rapidly propagate to the right of the initial pulse, and analysed these solutions to establish that the asymptotic solution obtained from the analysis of stationary KdV does not, in fact, exist. Instead, this solution must eventually decay on sufficiently long time scales (known as ``metastability'' \cite{Giardetti}). In this study, we instead consider the small-time behaviour of the system, focusing on the emergence of these oscillations rather than their long-time effects.

More general forms of the fifth-order KdV equation have been obtained by keeping higher-order terms related to nonlinear effects. These equations have been used to study solitary wave dynamics and wave radiation for free-surface flows 
\cite{benney1966long,grimshaw2002higher,horikis2021higher,hunter1983solitary,marchant1990extended,vanden1991elevation}, as well as internal layers \cite{grimshaw2002higher, khusnutdinova2018soliton} and ice sheets \cite{il2015soliton,marchenko1988long}. Other studies have studied dispersive shock waves \cite{sprenger2017shock} and examined the decay or persistence of soliton solutions \cite{tan2002semi} in fifth-order KdV systems.

Exponential asymptotic studies such as \cite{Grimshaw,grimshaw2010exponential,trinh2010exponential} examined the role of Stokes' phenomenon in asymptotic solutions to \eqref{e:kdv50}, showing that the exponentially small constant-amplitude oscillations are generated across curves in the analytically continued complex plane known as Stokes lines. These studies identified a one-parameter family of admissible wave profiles: while exponential asymptotics determines the jump in amplitude across a Stokes line, an additional condition is required to determine the amplitude on either side. Common choices include requiring the solution to be waveless on one side of the pulse or symmetric about its centre. Collectively, these works established that every solitary-wave solution of \eqref{e:kdv50} must possess oscillatory tails on one or both sides of the pulse. The work in \cite{Grimshaw} used energy arguments to show that a solution can persist indefinitely only when these oscillations are symmetric. This was supported by numerical and asymptotic work such as \cite{benilov1993generation}, which demonstrated that the waves naturally propagate on the right-hand side of the wave, and the pulse must decay in the large time limit.

\begin{figure}
\centering
\includegraphics[width=0.99\textwidth]{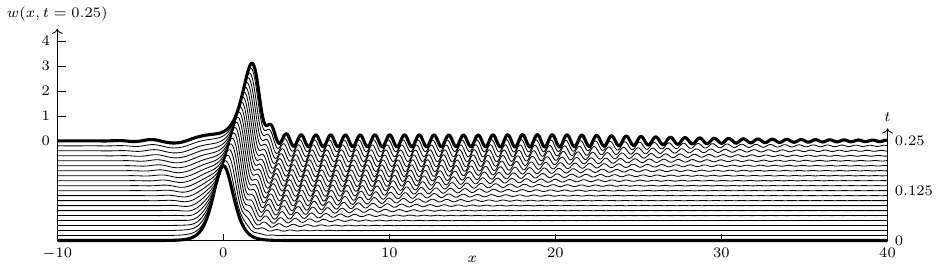}
\caption{Waterfall plot of the numerical solution of KdV5 with $c = 6$ and $\eps = 0.15$, for a range of $t$ from $t = 0$ to $t = 0.25$. The initial and final states are shown as thick black curves. The scale on the left shows the value of $w(x,t)$ at $t = 0.25$; this vertical scale is identical for all solutions. The scale on the right indexes the timesteps in the waterfall plot. The rapid propagation of waves to the right can be seen in these simulations.}
\label{fig:KdV5_num}
\end{figure}

These results naturally raise the question of how the generalised solitary wave emerges at small time. Numerical solutions of \eqref{e:kdv50} with initial data given by \eqref{e:kdvsoliton} show that exponentially small oscillations rapidly propagate to the right of the pulse \cite{benilov1993generation}, as shown in Figure \ref{fig:KdV5_num}. Further simulations indicate that the propagation speed depends on the small parameter $\eps$. The aim of this paper is to determine the mechanism governing this propagation using exponential asymptotic methods. Unlike the stationary problem, this requires an exponential asymptotic analysis of the full nonlinear partial differential equation \eqref{e:kdv50}, rather than the ordinary differential equations considered in previous studies.

Exponential asymptotic methods have previously been applied to linear partial differential equations (PDEs) \cite{chapman2005exponential,howls2004higher,body2005exponential}, revealing the central role played by Stokes' phenomenon in determining solution behaviour. These analyses showed that a complete asymptotic description requires tracking interactions between different exponential contributions, including second-generation Stokes' phenomenon, in which one exponentially small contribution switches on another, and higher-order Stokes' phenomenon, which can lead to Stokes lines terminating at points known as Stokes crossing points. These ideas have been applied to linear water wave problems, such as in \cite{lustri2013steady,lustri2022exponential}.

In this study, we adapt the approach of \cite{chapman2005exponential} for analysing exponentially small effects in PDEs. This method builds on the optimal truncation and matched asymptotic framework developed for linear differential equations in \cite{Daalhuis}, and its factorial-over-power extension to nonlinear differential equations in \cite{Chapman}. The central idea is to determine the behaviour of the $n$th term of an asymptotic series as $n\to\infty$, known as the late-order behaviour, and use this to truncate the series optimally. As shown by Berry \cite{berry1988stokes}, optimal truncation leaves an exponentially small remainder that can be studied directly. The methods of \cite{Daalhuis} can then be applied to this remainder to determine the exponentially small changes in behaviour that occur across Stokes lines.

Chapman \emph{et al.} \cite{chapman2005exponential} demonstrated that this framework extends naturally to the more complicated Stokes interactions that arise in PDEs. Second-generation Stokes switching can be captured through the late-order behaviour of the remainder series, while higher-order Stokes' phenomenon arises through Stokes switching in the prefactor of the late-order ansatz. These ideas were further developed in \cite{shelton2025exponential}, which provides a systematic framework for computing such effects.

Most existing applications of these techniques have focused on linear PDEs. Comparatively little work has been carried out for nonlinear PDEs \cite{chapman2007shock, mccue2026small}, largely because of the complexity of their Stokes structure. For example, \cite{chapman2007shock} examined the viscous Burgers equation, with viscosity as the small asymptotic parameter, and found that even this exactly solvable problem exhibits an intricate Stokes structure. \cite{mccue2026small} studied the KdV equation in the small-time limit and found that, although the real solution appears relatively simple, the analytically continued solution rapidly develops infinitely many singularities in the complex plane. Even numerical studies of Stokes surfaces for nonlinear PDEs, as in \cite{johnson2024three}, have proven challenging. These examples illustrate that nonlinear PDEs can generate highly complicated Stokes interactions.

In the present study, we navigate this complicated Stokes structure by constructing an auxiliary linear PDE that reproduces the relevant exponential asymptotic structure of \eqref{e:kdv50}. Specifically, we choose the auxiliary problem so that its late-order terms possess the same singulants and prefactors as those of the nonlinear equation. Consequently, the auxiliary problem reproduces the portion of the Stokes structure that controls the stationary oscillations studied here. The principal advantage of this approach is that powerful analytical tools, such as the method of steepest descents, can be applied directly to the linear problem but not easily to the nonlinear PDE. We use these techniques to identify the Stokes and higher-order Stokes lines that govern the solution behaviour of \eqref{e:kdv50}.

Having established the Stokes structure of the auxiliary problem, we then apply the factorial-over-power framework of \cite{chapman2005exponential} to both the linear and nonlinear equations. This allows us to determine an asymptotic description of the stationary waves and the region in which they are present, and hence their speed of propagation. We find that the nonlinear analysis differs only slightly from the corresponding linear calculation, so that much of the required information can be obtained from the preceding linear analysis.

The most significant difficulty in both problems is determining the multiplicative factor appearing in the prefactor of the late-order ansatz. This factor is usually obtained by matching the late-order behaviour with an inner solution near a complex singularity. For linear ordinary differential equations, the matching constant can often be determined exactly \cite{Daalhuis}, while nonlinear ordinary differential equations typically require a numerical calculation \cite{Chapman}. Even in more complicated settings, such as \cite{crew2024resurgent}, systematic analytical and numerical procedures exist for determining these quantities and the associated solution behaviour.

For partial differential equations, however, the situation is considerably more complicated because the multiplicative factor is no longer a constant but a function of a characteristic variable. Furthermore, numerical evaluation of this function is not, by itself, sufficient to determine the relevant solution behaviour. Second-generation Stokes switching is generated by singularities of the multiplicative factor, meaning that the crucial information is not its global form but rather the locations and local behaviour of these singularities. Extracting this information directly from numerical calculations is challenging. Instead, we exploit the requirement that the solution is initially waveless, which confines the oscillations to a finite region of the solution domain. This condition alone is sufficient to determine the stationary-wave behaviour.

\section{Formulating the Problem}

\subsection{Fifth-Order KdV Equation (KdV5)}
We consider the singularly-perturbed fifth-order KdV equation, given by \eqref{e:kdv50}, with the initial condition
\begin{equation}\label{e:ic0}
w(x,0) = \tfrac{c}{2}\sech^2\left(\tfrac{\sqrt{c}}{2}x\right).
\end{equation}
We will also apply the condition that $w(x,t) \to 0$ in the limit $x \to \pm\infty$. In the unperturbed KdV equation \eqref{e:kdv} this initial condition produces the one-soliton solution \eqref{e:kdvsoliton}. Motivated by this, we will move into a travelling wave frame by setting $\xi = x - c t$, but we will retain the time variable $t$ so that $u \equiv u(\xi,t)$, giving
\begin{align}\label{e:kdv5}
\eps^2 w_{\xi\xi\xi\xi\xi} + w_{\xi\xi\xi}  + 6 ww_{\xi} - c w_{\xi} + w_t = 0,& \qquad (\xi,t) \in (-\infty,\infty) \times [0,\infty),\\
w(\xi,0) = \tfrac{c}{2}\sech^2\left(\tfrac{\sqrt{c}}{2}\xi\right),& \qquad \xi \in (-\infty,\infty),\label{e:ic}\\
w(\xi,t) \to 0 &\qquad \mathrm{as} \quad \xi \to \pm\infty.\label{e:kdv5-ff}
\end{align}
Expanding the solution as a power series in $\eps$ gives
\begin{equation}\label{e:series}
w(\xi,t) \sim \sum_{n=0}^{\infty}\eps^{2n}w^{(n)}(\xi,t)\quad \as \quad \eps \to 0.
\end{equation}
The leading-order solution  $w^{(0)}$ satisfies the KdV equation \eqref{e:kdv}, and the initial condition \eqref{e:ic} therefore produces the soliton solution \eqref{e:kdvsoliton}, so that
\begin{equation}\label{e:u0}
w^{(0)}(\xi,t) = \tfrac{c}{2}\sech^2\left(\tfrac{\sqrt{c}}{2}\xi\right).
\end{equation}
We will refer to the localised solitary wave associated with this leading-order solution as the ``primary solitary wave''. The recurrence relation for $w^{(n)}$ is obtained by substituting \eqref{e:series} into \eqref{e:kdv5} and balancing terms at successive powers of $\eps$. For $n\geq 1$, we find
\begin{align}\label{e:kdv5recur}
w^{(n-1)}_{\xi\xi\xi\xi\xi} + w^{(n)}_{\xi\xi\xi} + 6 \sum_{j=0}^{n}w^{(j)} w^{(n-j)}_{\xi}-c w^{(n)}_{\xi}+  w^{(n)}_{t} = 0, & \qquad (\xi,t) \in (-\infty,\infty) \times [0,\infty),\\
w^{(n)}(\xi,0) = 0,& \qquad \xi \in (-\infty,\infty),\label{eq:kdv5recur-ic}\\
w^{(n)}(\xi,t) \to 0& \qquad \mathrm{as} \quad \xi \to \pm\infty.\label{eq:kdv5recur-farfield}
\end{align}
Typically, applying an exponential asymptotic analysis would involve determining the large-$n$ limit of the series terms $w^{(n)}$; however, applying this approach naively does not work in the current setting. This analysis would simply reproduce the result from \cite{Grimshaw,trinh2010exponential} with no additional exponential contributions. Such a result cannot possibly satisfy the boundary condition as $\xi \to \pm\infty$, as it contains exponentially small waves of constant amplitude on one or both sides of the primary solitary wave.

\subsection{Auxiliary Linear PDE (lKdV5)}

The fact that \eqref{e:kdv5} is a nonlinear PDE means that it is challenging to study directly using exponential asymptotic methods. Instead of proceeding directly to an exponential asymptotic analysis of \eqref{e:kdv5} we first formulate an auxiliary linear problem with carefully chosen properties. We will therefore consider 
\begin{align}\label{e:kdv5lin}
\eps^2 u_{\xi\xi\xi\xi\xi} + u_{\xi\xi\xi} + u_t = \left[\tfrac{c}{2}\sech^2\left(\tfrac{\sqrt{c}}{2}\xi\right)\right]_{\xi\xi\xi},& \qquad (\xi,t) \in (-\infty,\infty) \times [0,\infty),\\
u(\xi,0) = \tfrac{c}{2}\sech^2\left(\tfrac{\sqrt{c}}{2}\xi\right),& \qquad \xi \in (-\infty,\infty),\\
u(\xi,t) \to 0& \qquad \mathrm{as} \quad \xi \to \pm\infty.\label{e:kdv5lin-ff}
\end{align}
This linear equation is selected because its asymptotic solutions possess the same Stokes structure as those of \eqref{e:kdv5} in the relevant regime. We will show this claim to be true in Section \ref{s:KdV5}. 

At this stage, one may naively expand the solution to lKdV5 as a series as in \eqref{e:series} and find
\begin{equation}\label{e:u0l}
u^{(0)}(\xi,t) = \tfrac{c}{2}\sech^2\left(\tfrac{\sqrt{c}}{2}\xi\right).
\end{equation}
It is possible to obtain a recurrence relation similar to \eqref{e:kdv5recur}--\eqref{eq:kdv5recur-farfield} and calculate further correction terms, and even to perform an exponential asymptotic analysis to study Stokes' phenomenon in this system. However, as in KdV5, an exponential asymptotic analysis on this system will not capture time-varying exponentially small behaviour. This suggests that the solution contains a small-time boundary layer, in which the overall system rapidly approaches a known steady solution. 

\subsection{Small-time scaling}

To capture the rapid spread of ripples in the solution to KdV5 \eqref{e:kdv5}--\eqref{e:kdv5-ff} and lKdV5 \eqref{e:kdv5lin}--\eqref{e:kdv5lin-ff}, we must apply an appropriate small-time scaling to the system. This scaling should be chosen such that the term in each equation containing a derivative with respect to time contributes to the exponential asymptotic analysis. This occurs if $t = \mathcal{O}(\eps^2)$ in the limit $\eps \to 0$, as confirmed by the later analysis in Sections \ref{s:lKdV5} and \ref{s:KdV5}.

This scaling could alternatively be determined from the exponential contributions calculated in \cite{mccue2026small} for the unperturbed KdV equation with small time. Balancing these exponentials with their fifth derivative (scaled by $\eps^4$) gives the desired $\mathcal{O}(\eps^2)$ scaling for $t$. 

\subsubsection{KdV5}

Applying the rescaling $t = \eps^2 \tau$ to \eqref{e:kdv5}--\eqref{e:kdv5-ff} and letting $w(\xi,t) \equiv u(\xi,\tau)$ gives
\begin{align}\label{e:kdv5ST}
\eps^4 w_{\xi\xi\xi\xi\xi} + \eps^2\left( w_{\xi\xi\xi}  + 6 w w_{\xi} - c w_{\xi} \right) +w_{\tau} = 0, & \qquad (\xi,\tau) \in (-\infty,\infty) \times [0,\infty),\\
w(\xi,0) = \tfrac{c}{2}\sech^2\left(\tfrac{\sqrt{c}}{2}\xi\right),& \qquad \xi \in (-\infty,\infty),\label{e:kdv5ST-ic}\\
w(\xi,\tau) \to 0 &\qquad \mathrm{as} \quad \xi \to \pm\infty.\label{e:kdv5ST-ff}
\end{align}
Expanding the solution as a power series, as in \eqref{e:series}, gives 
\begin{equation}
w^{(0)}(\xi,\tau) =  \tfrac{c}{2}\sech^2\left(\tfrac{\sqrt{c}}{2}\xi\right).
\end{equation}
Note that this solution is steady in $\tau$. This highlights that the small-time boundary layer has not been selected to balance the $u_{\tau}$ term with the leading-order contribution, as is typical in boundary layer problems, but rather to balance the exponentially small effects. For $n \geq 1$, the recurrence relation becomes (with the convention that $u^{(-1)} = 0$)
\begin{align}\label{e:kdv5recurST}
w^{(n-2)}_{\xi\xi\xi\xi\xi} + w^{(n-1)}_{\xi\xi\xi} + 6 \sum_{j=0}^{n-1}w^{(j)} w^{(n-1-j)}_{\xi}-c w^{(n-1)}_{\xi}+  w^{(n)}_{\tau} = 0,& \qquad (\xi,\tau) \in (-\infty,\infty) \times [0,\infty),\\
w^{(n)}(\xi,0) = 0,& \qquad \xi \in (-\infty,\infty),\label{eq:kdv5recurST-ic}\\
w^{(n)}(\xi,\tau) \to 0& \qquad \mathrm{as} \quad \xi \to \pm\infty.\label{eq:kdv5recurST-farfield}
\end{align}
By direct integration, we find that $w^{(1)}(\xi,\tau) = 0$ and
\begin{align}\label{eq:u2kdv5}
w^{(2)}(\xi,\tau) = &\tfrac{c^{7/2}}{4}\tau\left(17\sech^6\big(\tfrac{\sqrt{c}}{2}\xi\big)\tanh\big(\tfrac{\sqrt{c}}{2}\xi\big)- 26\sech^4\big(\tfrac{\sqrt{c}}{2}\xi\big)\tanh^3\big(\tfrac{\sqrt{c}}{2}\xi\big) + 2 \sech^2\big(\tfrac{\sqrt{c}}{2}\xi\big)\tanh^5\big(\tfrac{\sqrt{c}}{2}\xi\big)\right).
\end{align}
Subsequent terms can be obtained by repeatedly applying \eqref{e:kdv5recurST}--\eqref{eq:kdv5recurST-farfield}. 

\subsubsection{lKdV5}\label{ss:lkdv5}

Applying the rescaling $t = \eps^2 \tau$ to \eqref{e:kdv5lin}--\eqref{e:kdv5lin-ff} and letting $w(\xi,t) \equiv u(\xi,\tau)$ gives
\begin{align}\label{e:lkdv5}
\eps^4 u_{\xi\xi\xi\xi\xi} + \eps^2 u_{\xi\xi\xi} + u_{\tau} = \eps^2 \left[  \tfrac{c}{2}\sech^2\left(\tfrac{\sqrt{c}}{2}\xi\right)\right]_{\xi\xi\xi},& \qquad (\xi,\tau) \in (-\infty,\infty) \times [0,\infty),\\
u(\xi,0) = \tfrac{c}{2}\sech^2\left(\tfrac{\sqrt{c}}{2}\xi\right),& \qquad \xi \in (-\infty,\infty),\label{e:lkdv5-ic}\\
u(\xi,\tau) \to 0 &\qquad \mathrm{as} \quad \xi \to \pm\infty.\label{e:lkdv5-ff}
\end{align}
Expanding the solution as a power series, as in \eqref{e:series}, gives 
\begin{equation}
u^{(0)}(\xi,\tau) =  \tfrac{c}{2}\sech^2\left(\tfrac{\sqrt{c}}{2}\xi\right),
\end{equation}
as is the case for KdV5. For $n \geq 1$, the recurrence relation becomes
\begin{align}\label{e:lkdv5recur}
u^{(n-2)}_{\xi\xi\xi\xi\xi} + u^{(n-1)}_{\xi\xi\xi}+  u^{(n)}_{\tau} = 0,& \qquad (\xi,\tau) \in (-\infty,\infty) \times [0,\infty),\\
u^{(n)}(\xi,0) = 0,& \qquad \xi \in (-\infty,\infty),\label{e:lkdv5recur-ic}\\
u^{(n)}(\xi,\tau) \to 0& \qquad \mathrm{as} \quad \xi \to \pm\infty.\label{e:lkdv5recur-farfield}
\end{align}
By direct integration, we find that the first correction term does not cancel as before, and is instead
\begin{align}\label{eq:u1l}
u^{(1)}(\xi,\tau) = -\tfrac{c^{5/2}}{2}\tau\left(2\sech^4\big(\tfrac{\sqrt{c}}{2}\xi\big)\tanh\big(\tfrac{\sqrt{c}}{2}\xi\big)- \sech^2\big(\tfrac{\sqrt{c}}{2}\xi\big)\tanh^3\big(\tfrac{\sqrt{c}}{2}\xi\big)\right).
\end{align}
Subsequent terms can be obtained by repeatedly applying \eqref{e:kdv5recur}--\eqref{eq:kdv5recur-farfield}.

\section{lKdV5: Steepest Descents}\label{s:linSD}

As lKdV5 is a linear PDE, we may study it using the method of steepest descents. We will take the Fourier transform with respect to $\xi$, and use the pair
\begin{equation}\label{e:fourier}
\hat{f}(k) = \int_{-\infty}^{\infty}f(\xi)\e^{-\i k \xi}\d \xi,\qquad f(x) = \frac{1}{2\pi}\int_{-\infty}^{\infty}\hat{f}(\xi)\e^{\i k \xi}\d k.
\end{equation}
Taking the Fourier transform of \eqref{e:lkdv5}--\eqref{e:lkdv5-ff} results in
\begin{equation}\label{e:kdvlinft}
\i\eps^2 k^5 \hat{u} - \i k^3 \hat{u} + \diff{\hat{u}}{t} = -\i k^3 \hat{u}_0, \qquad \hat{u}(0) = \hat{u}_0.
\end{equation}
Solving this differential equation gives 
\begin{equation}\label{e:kdvlinuhat}
\hat{u} = \left(\eps^2 k^2 \e^{-\i k^3 (\eps^2 k^2 - 1)t} - 1\right)\frac{\hat{u}_0}{\eps^2 k^2 - 1},\qquad \hat{u}_0 = \frac{2\pi k} {\sinh\left(k\pi/\sqrt{c}\right)}.
\end{equation}
Inverting \eqref{e:kdvlinuhat} gives 
\begin{equation}\label{e:kdvlinint}
u = \frac{1}{2\pi}\int_{-\infty}^{\infty} \left(\eps^2 k^2 \e^{-\i k^3 (\eps^2 k^2 - 1)t} - 1\right)\frac{\hat{u}_0 \e^{\i k \xi}}{\eps^2 k^2 - 1} \d k,
\end{equation}
where we indent the integral below the poles at $k = \pm 1/\eps$, which forces the oscillations to appear for positive $\xi$. Splitting the integral between the steady and the unsteady contribution gives $u = v + u_s$, where
\begin{align}
\label{e:kdvlinv0}v &= \frac{1}{2\pi}\int_{-\infty}^{\infty}  \frac{\eps^2 k^2 \hat{u}_0\e^{-\i k^3 (\eps^2 k^2 - 1)t + \i k\xi}}{\eps^2 k^2 - 1}   \d k,\\
\label{e:kdvlinus}u_s &= -\frac{1}{2\pi}\int_{-\infty}^{\infty} \frac{\hat{u}_0\e^{\i k \xi} }{\eps^2 k^2 - 1} \d k.
\end{align}

\subsection{Steady Contribution}
The asymptotic behaviour of \eqref{e:kdvlinus} is straightforward, and similar integrals have been calculated in numerous previous studies. Using the expression for $\hat{u}_0$ from \eqref{e:kdvlinuhat}, 
the transformation $k = s/\eps$, and exploiting the symmetry of the integrand gives
\begin{equation}
\label{e:kdvlinussplit}u_s =-\frac{2}{\eps^2}\int_{0}^{\infty}  \frac{s \e^{\i s (\xi + \i\pi/\sqrt{c})/\eps}}{(s^2 - 1)( 1+ \e^{-2s\pi/(\sqrt{c}\eps)})}  \d s + \mathrm{c.c.},
\end{equation}
where c. c. refers to the complex conjugate contribution. The asymptotic behaviour of this integral contains a contribution from the end-point at $s = 0$, which contributes for all values of $\xi$, and the pole at $s = 1$, which contributes if $\xi > 0$. The end-point contribution is a power series in $\eps$, which correspond to localised contributions from the soliton and its higher-order corrections.
The pole contribution is 
\begin{equation}\label{eq:steadywaves}
u_{s,p} \sim -\frac{2\pi\i}{\eps^2}  \e^{\i(\xi + \i\pi/\sqrt{c})/\eps} + \mathrm{c. c.} = \frac{4\pi}{\eps^2}\e^{-\pi/(\sqrt{c}\eps)}\sin\left(\frac{\xi}{\eps}\right).
\end{equation}
This gives oscillations for $\xi > 0$, as expected.

\subsection{Unsteady Contribution}

Using the expression for $\hat{u}_0$ from \eqref{e:kdvlinuhat}, the unsteady integral \eqref{e:kdvlinv0} becomes
\begin{equation}
\label{e:kdvlinv}v = \int_{-\infty}^{\infty}  \frac{\eps^2 k^3 \e^{-\i k^3 (\eps^2 k^2 - 1)t + \i k \xi}}{(\eps^2 k^2 - 1)(\e^{k\pi/\sqrt{c}} + \e^{-k\pi/\sqrt{c}})}  \d k
\end{equation}
Setting $k = s/\eps$ and $t = \eps^2 \tau$ gives $
v 
= I(\xi,\tau) + \mathrm{c. c.},
$ 
where
\begin{equation}\label{e:kdvlinI}
I(\xi,\tau) = \frac{2}{\eps^2}\int_{0}^{\infty}  \frac{s^3  \e^{(-\i s^3 (s^2 - 1)\tau + \i s (\xi +\i \pi/\sqrt{c}))/\eps}}{(s^2 - 1)(1 + \e^{-2s\pi/(\sqrt{c}\eps)})}  \d s,
\end{equation}
This introduces the small time scaling that we expect to see in the solution. We can now apply the usual method of steepest descents. By setting 
\begin{equation}\label{e:gs}
g(s) = -\i s^3 (s^2 - 1)\tau + \i s (\xi +\i \pi/\sqrt{c} ),
\end{equation} 
and solving $g'(s) = 0$, we identify that the integrand contains four saddle points. We are primarily concerned with the saddle points in the right half plane, located at $s = s_1$ and $s = s_2$ where
\begin{equation}
s_{1} = \left(\tfrac{3}{\sqrt{10}} + \tfrac{1}{\sqrt{10}}\sqrt{9 + 20\left(\xi + \tfrac{\i\pi}{\sqrt{c}}\right)\tfrac{1}{\tau}}\right)^{1/2},\qquad s_{2} = \left(\tfrac{3}{\sqrt{10}} - \tfrac{1}{\sqrt{10}}\sqrt{9 + 20\left(\xi + \tfrac{\i\pi}{\sqrt{c}}\right)\tfrac{1}{\tau}}\right)^{1/2}.
\end{equation}
$I(\xi,\tau)$ contains four possible asymptotic contributions: one from the end-point at $s = 0$; one from the each saddle point at $s = s_1$ and $s=s_2$, and one from the pole at $s = 1$. We define a useful set of variables,
\begin{equation}\label{e:newvars}
\sigma = \xi + \frac{\i\pi}{\sqrt{c}}, \qquad  \nu = \frac{\sigma}{\tau}, \qquad \alpha = 3 + \sqrt{9 + 20\nu}, \qquad \tilde{\alpha} = 3 - \sqrt{9 + 20\nu}.
\end{equation} 
In terms of these variables, we have
\begin{equation}
s_{1} = \left({\tfrac{1}{\sqrt{10}}\alpha}\right)^{1/2},\qquad s_{2} = \left({\tfrac{1}{\sqrt{10}}\tilde{\alpha}}\right)^{1/2}.
\end{equation}

\begin{figure}
\centering
\subfloat[$\xi = -1$, $\tau = 1$]{
\includegraphics[width=0.3\textwidth]{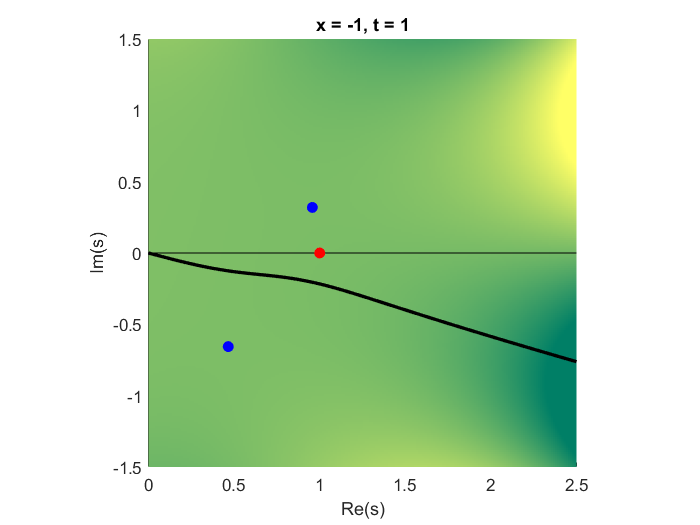}
}
\subfloat[$\xi = 0.25$, $\tau = 1$]{
\includegraphics[width=0.3\textwidth]{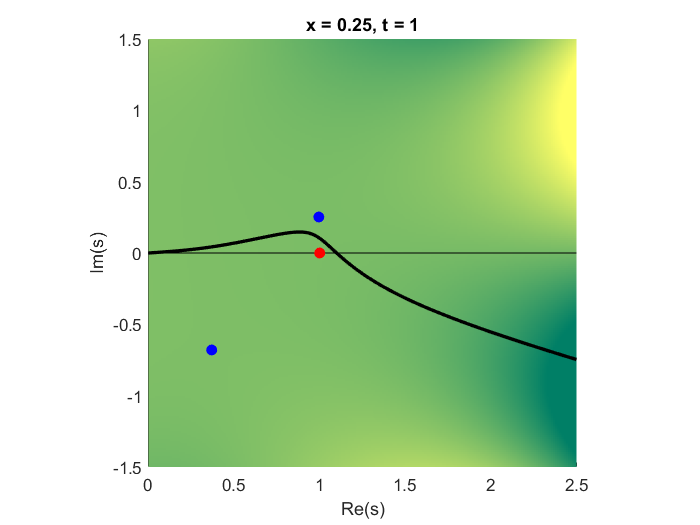}
}
\subfloat[$\xi = 0.75$, $\tau = 1$]{
\includegraphics[width=0.3\textwidth]{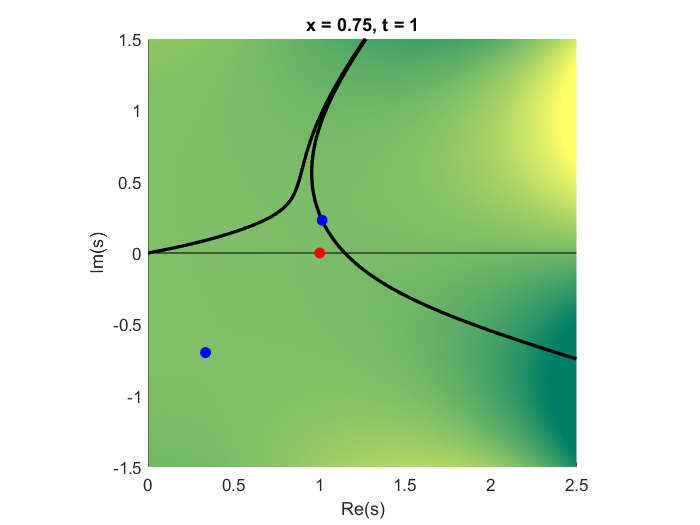}
}

\subfloat[$\xi = 0.75$, $\tau = 3$]{
\includegraphics[width=0.3\textwidth]{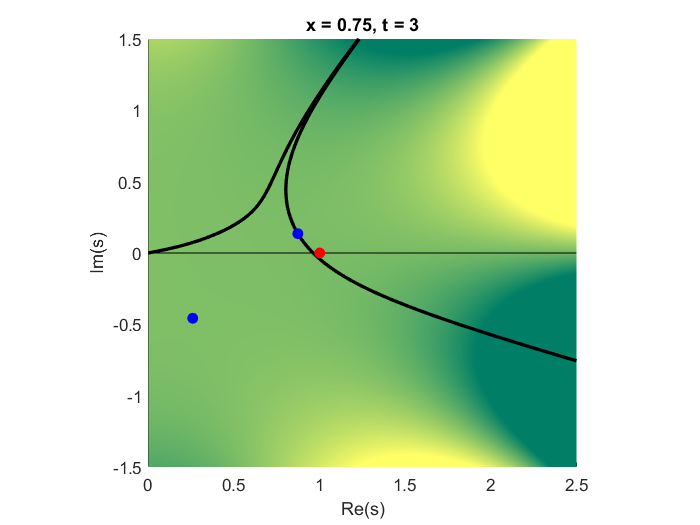}
}
\subfloat[$\xi = -1$, $\tau = 3$]{
\includegraphics[width=0.3\textwidth]{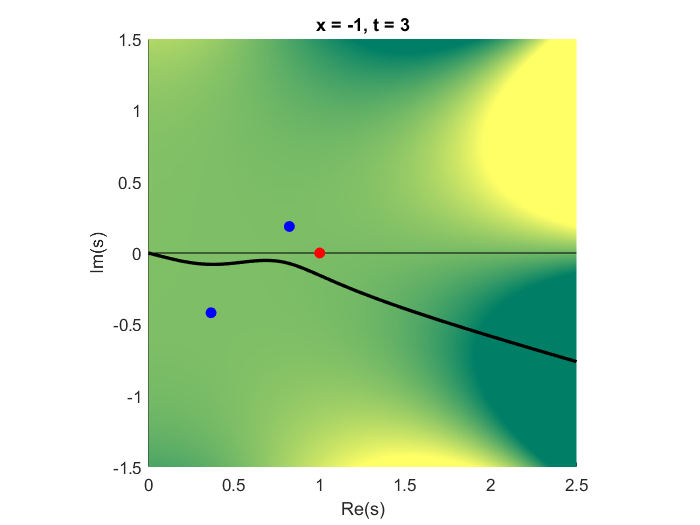}
}
\subfloat[$\xi = -4$, $\tau = 1$]{
\includegraphics[width=0.3\textwidth]{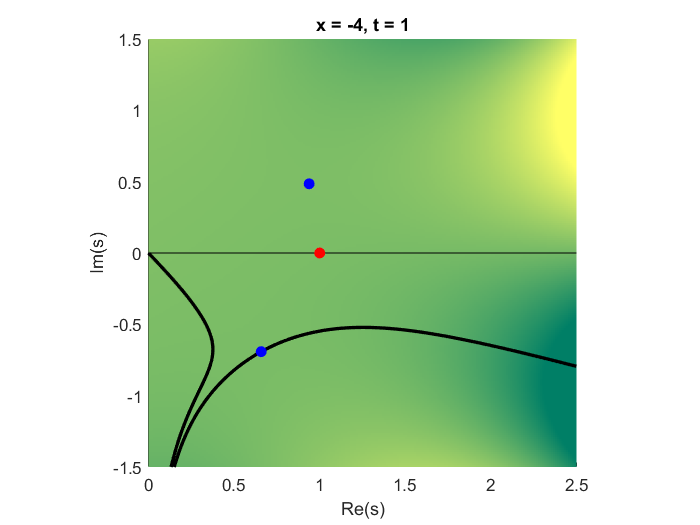}
}
\caption{Steepest descent curves for integral $I(\xi,\tau)$ from \eqref{e:kdvlinI} with $c = 1$, shown over $\mathrm{Re}(g(s))$. Saddle points are shown as blue circles, and the pole is shown as a red circle. The steepest descent curve is a solid black line that begins at $s = 0$. Darker regions have more negative values of $\mathrm{Re}(g(s))$. In (a) the only asymptotic contribution comes from the end-point. In (b) the end-point has switched on the pole contribution. In (c) the end-point has switched on the upper-half saddle-point contribution, so the solution contains end-point, pole, and saddle-point contributions. In (d) the saddle point has switched off the pole contribution, leaving only saddle-point and end-point contributions. In (e) the end point has switched off the saddle-point contribution, so only the end-point contribution remains. In (f) the end point has switched on the lower-half saddle contribution. }
\label{fig:v_SD_curves}
\end{figure}

A steepest descent analysis can identify which contributions are present in each region of the ($\xi$, $\tau$) plane. Samples of the different possible configurations are presented in Figure \ref{fig:v_SD_curves}. Figures (a) and (b) depict steepest descent contours on either side of the pole, so the contribution in (a) does not pick up the pole contribution but the contour in (b) does. The contours in (c) and (d) both pass through the saddle point at $s_1$; the contour in (c) picks up the pole contribution, but the contour in (d) does not. The contour in (f) shows that the steepest descent curves can also pass through the saddle point at $s_2$, although this will not be important in determining how the exponentially small waves behave, as the steepest curve never passes through $s_2$ while also including the pole contribution.

\subsubsection{Stokes Structure}

\begin{figure}
\centering
\includegraphics{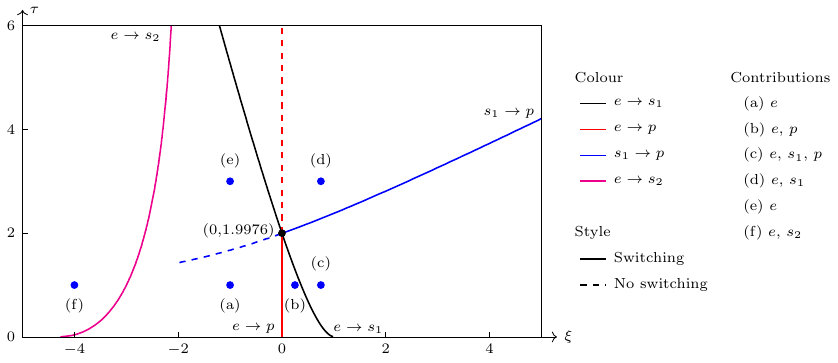}
\caption{Stokes' phenomenon in the solution to lKdV5 with $c = 1$. The solid lines represent Stokes lines across which a dominant contribution switches on or off a subdominant contribution. There are four possible contributions represented: the end-point ($e$), the two saddle points ($s_1$, $s_2$), and the pole ($p$). This is indicated in the legend by an arrow; for example, $e \to p$ indicates a Stokes curve where the pole contribution is switched on or off by the end-point contribution. Dashed lines represent curves where the two contributions have equal phase, but no switching occurs. This is a consequence of the Stokes lines being non-adjacent or the dominant contribution not being present in the solution; this indicates the presence of higher-order Stokes' phenomenon. The six points marked by a blue circle correspond to the sample steepest descent curves in Figure \ref{fig:v_SD_curves}. The blue dashed line does not end but continues to the left. The remaining portion of the line was omitted for visual clarity.}
\label{fig:v_stokes}
\end{figure}

We can determine the full Stokes structure, and hence the asymptotic contributions for different values of $\xi$ and $\tau$, using the Stokes conditions. A contribution from $s = a$ can switch a contribution from $s = b$ if 
\begin{equation}\label{e:stokescond}
\Im\left[g(a)\right] = \Im\left[g(b)\right], \qquad \Re\left[g(a)\right] > \Re\left[g(b)\right].
\end{equation}
We can use this condition to identify any possible Stokes lines in the asymptotic behaviour as $(\xi,\tau)$ vary. Note that these conditions alone are not sufficient to guarantee Stokes' phenomenon, but the steepest descent contours can easily be checked to determine whether any Stokes switching occurs across curves that satisfy \eqref{e:stokescond}. This produces the Stokes structure depicted in Figure \ref{fig:v_stokes}.

This figure shows which terms contribute in the six distinct regions, denoted (a)--(f). These correspond to sample steepest descent curves in Figure \ref{fig:v_SD_curves}. The system contains a Stokes line across which the end-point switches each of the two saddles. It also contains two additional Stokes lines, across which the pole is switched by the end-point and the saddle. 

\subsubsection{Higher-Order Stokes' Phenomenon}

\begin{figure}
\centering
\includegraphics{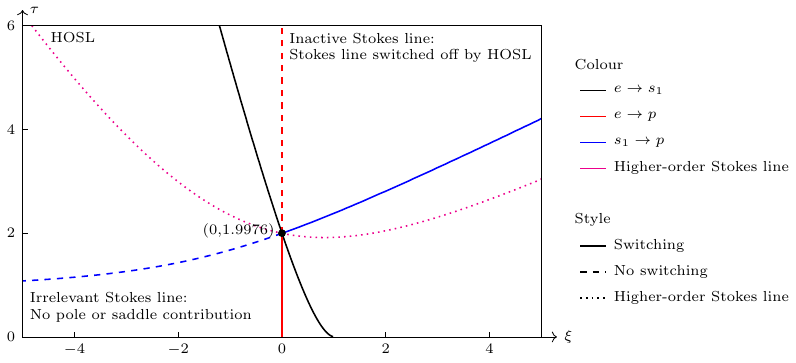}
\caption{The higher-order Stokes line (HOSL) in the solution to the lKdV5 equation with $c = 1$ due to the colinearity condition between the end-point, saddle and pole contributions. Below this curve, the saddle-point/pole Stokes line is irrelevant, as the saddle-point and pole contributions are inactive, so no Stokes switching can occur. Above this curve, the end-point/pole Stokes line is inactive, as the two contributions are not adjacent in the steepest descent curve. This is required for the solution to be continued around the Stokes crossing point at $(\xi,\tau) \approx (0,1.9976)$.}
\label{fig:v_HOSL}
\end{figure}

Figure \ref{fig:v_stokes} contains a curve that satisfies \eqref{e:stokescond} but across which no switching occurs, represented by a dashed blue line. The curve satisfying
\begin{equation}\label{e:SLsd}
\mathrm{Im}[g(s_1)] = \mathrm{Im}[g(1)]
\end{equation}
does not produce Stokes switching if $\xi < 0$, even though $\mathrm{Re}[g(s_1)] > \mathrm{Re}[g(1)]$. This occurs because the saddle contribution is not active in the solution, so the Stokes line is inactive. Hence there must be a higher-order Stokes line in the solution \cite{chapman2005exponential,howls2004higher} that intersects $(\xi,\tau) \approx(0,1.9976)$, across which the Stokes line becomes inactive. The solution contains a higher-order Stokes line on the curve which satisfies the 
colinearity condition
\begin{equation}\label{e:HOSLsd}
\mathrm{Im}\left[\frac{g(s_1) - g(0)}{g(1) - g(0)}\right] = 0,
\end{equation}
which simplifies to give $\mathrm{Im}\left[{g(s_1)}/{g(1)}\right] = 0$. 
Below this higher-order Stokes line, 
the
Stokes line across which the end-point contribution switches on the pole contribution behaves as normal. Above this higher-order Stokes line,
no switching occurs
even though the usual Stokes line switching conditions are met. The higher-order Stokes line is shown in Figure \ref{fig:v_HOSL}.

\subsubsection{Asymptotic contributions}

Each of the three asymptotic contributions can be calculated using standard asymptotic integral methods (eg. from \cite{bender2013advanced}). The relevant calculations are summarised below, with all limits taken to be $\eps \to 0$:

\begin{itemize}
\item The end-point $I_e$ generates an asymptotic power series that reproduces the contribution developed in Section \ref{ss:lkdv5}.
\item The contribution from the pole $I_p$ may be found using the residue theorem, giving
\begin{equation}
I_p \sim \frac{2\pi\i}{\eps^2}\e^{\i(\xi + \i\pi/\sqrt{c})/\eps} + \mathrm{c. c.}.
\end{equation}
This precisely cancels the oscillations in the steady solution \eqref{eq:steadywaves}. 

\item The saddle contributions $I_{s,j}$ can be calculated using the saddle point formula \cite{bender2013advanced}, giving
\begin{equation}
I_{s_j} \sim \frac{2}{\eps^2}\left({-\frac{2\pi\eps}{g''(s_j)}}\right)^{1/2}\frac{s_j^3}{s_j^2 -1}\,\e^{g(s_j)/\eps} \quad \as\quad\eps \to 0,
\end{equation}
for $j = 1, 2$. It will be useful to calculate $I_{s_1}$ explicitly in order to study the Stokes switching associated with the saddle at $s_1$. In terms of the variables from \eqref{e:newvars}, we have
\begin{equation}\label{e:p(nu)}
g(s_1) = -p(\nu)\tau \quad \mathrm{where} \quad \qquad p(\nu) =-\frac{\i\alpha^{1/2}(20\nu + \alpha)}{25\sqrt{10}}.
\end{equation}
The new quantity $p(\nu)$ will be useful later. We also find
\begin{equation}
g''(s_1) = -\i\tau\alpha^{1/2}\sqrt{\tfrac{18}{5} + 8 \nu}.
\end{equation}
This means we have all the quantities required to explicitly calculate the saddle contributions. 
After some simplification, we find that
\begin{equation}\label{e:sdsaddle}
I_{s_1} \sim \frac{(-2)^{3/4}\alpha^{5/4}\sqrt{\pi}}{5^{1/4}\eps^{3/2}\tau^{1/2}(\alpha - 10)\sqrt{\alpha - 3}}\e^{-p(\nu)\tau/\eps} \quad \as \quad \eps \to 0.
\end{equation}
This is singular at $\nu = 2$, where the local behaviour is given by
\begin{equation}
I_{s_1} \sim \frac{2\e^{3\pi\i/4}\sqrt{7\pi}}{\eps^{3/2}\tau^{1/2}(\nu - 2)} \quad \as \quad \eps \to 0, \, \nu \to 2.
\end{equation}
It is worthwhile to understand the behaviour of this transient term in two limits: the limit that $\tau \to \infty$ for fixed $\xi$, and the limit that $\xi \to \infty$ for fixed $\tau$. If we fix $\xi$, we find that
\begin{equation}
I_{s_1} \sim -\frac{3^{3/4}}{5^{1/4}\eps^{3/2}}\sqrt{-\frac{\pi\i}{\tau}} \mathrm{exp}\left(\frac{ 3^{1/2}\i(6\tau/25 + (\xi + \i\pi/\sqrt{c}))}{5^{1/2}\eps}\right) \quad \mathrm{as} \quad \eps \to 0,\, \tau \to \infty.
\end{equation}
This implies that the saddle contribution is transient, and decays with $\mathcal{O}(\tau^{-1/2})$ as $\tau \to \infty$ for any fixed $\xi$. Hence, it does not persist outside of the small-time boundary layer, as we would expect.

To determine the far-field behaviour as $\xi \to \infty$ caused by the saddle, it is sufficient to show that
\begin{equation}
\mathrm{Re}[g(s_1)] \sim -\frac{\pi\xi^{1/4}}{5^{1/4}\sqrt{c}\tau^{1/4}} \quad \mathrm{as} \quad \xi \to \infty.
\end{equation}
This implies that the amplitude of the oscillations must decay as $\xi \to \infty$ for any fixed value of $\tau$, meaning that the solution does not contain growing or constant-amplitude waves in the far-field limit. The boundary condition that $u \to 0$ as $\xi \to \infty$ is therefore satisfied.

\end{itemize}

\subsection{Combined Solution}
\begin{figure}
\centering
\includegraphics{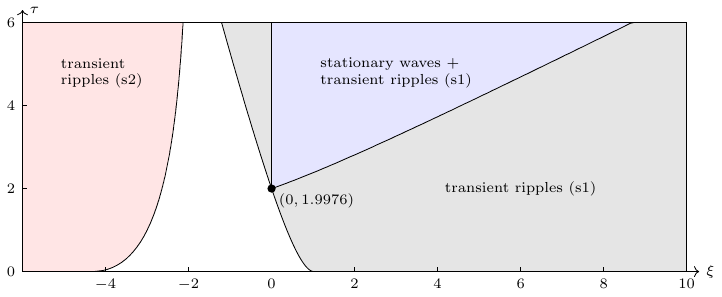}
\caption{Asymptotic contributions to $u(\xi,\tau)$ with $c = 1$. The white region contains no exponentially small contributions. The grey region represents the region containing decaying transient ripples from the upper-half saddle point. The red region represents the region containing decaying transient ripples from the lower-half saddle point. The blue region represents the region containing both transient ripples from the upper-half saddle point and stationary exponentially small oscillations that persist beyond the small-time boundary layer. This region appears at $\tau = 2$ at $\xi = 0$ and expands to the right such that the right edge of the ripple region follows $\xi \sim 2\tau$ as $\tau \to \infty$.}
\label{fig:ripple_region}
\end{figure}

By combining $v$ and $u_s$, we can identify the regions of the $(\xi,\tau)$ plane which contain oscillations. Neither of the two integrals possess pole contributions for $\xi < 0$, so we consider only the region $\xi > 0$.  If the pole contributions from both $v$ and $u_s$ are active, they cancel each other precisely. Hence, stationary oscillations are only present in the region in which $v$ does not have an integral contribution. This is shown in Figure \ref{fig:ripple_region}.

An important observation is that the initial condition can only be waveless if the constant-amplitude stationary oscillations are switched off in regions (b) and (c). This means that the contributions switched on by the end point/pole Stokes line and switched off by the saddle point/pole Stokes line must cancel precisely in order for the initial condition to be satisfied. This is evident in the steepest descent formulation, but will be a valuable observation for the later factorial-over-power analysis.

By direct calculation, it is possible to solve $\mathrm{Im}[g(s_1)] = \mathrm{Im}\mathrm{[g(1)]}$, or
\begin{equation}\label{e:StokesCond}
\mathrm{Im}\left[- \i s_1^3 (s_1^2 - 1)\tau + \i s_1 \left(\xi + \tfrac{\i \pi}{\sqrt{c}}\right) \right] = \mathrm{Im}\left[\i\left( \xi  + \tfrac{\i \pi}{\sqrt{c}}\right)\right]
\end{equation}
to determine the position of the right edge of the ripple region, which we denote $\xi = \xi_r$. For large $\xi$ and $\tau$ of comparable order, this eventually reduces to solving
\begin{equation}
\left(3\tau + \sqrt{\tau(9\tau + 20 \xi_r)}\right)^{3/2} \sim 25\xi_r \sqrt{10\tau} - 20 \xi_r \sqrt{3\tau + \sqrt{\tau(9\tau + 20\xi_r)}} \quad \as \quad \tau \to \infty.
\end{equation}
This can be solved to give $\xi_r \sim 2\tau$ as $\tau \to \infty$. Higher-order corrections to $\xi_r$ can be obtained, but we will not provide these calculations here. Hence, the position of the expanding front of ripples is $x \sim 2t/\eps^2$ in the limit that $\eps \to 0$ and $t \to \infty$. 

\subsection{Numerical Comparisons}

\begin{figure}
\centering
\subfloat[Numerical solution (blue): $c = 2$, $\eps = 0.125$, $t = 0.5$.]{
\includegraphics[width=0.9\textwidth]{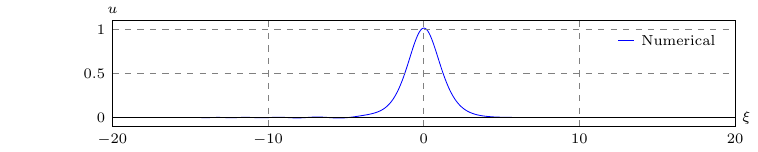}
 }
 
\subfloat[Numerical solution (blue): exponentially small oscillations]{
\includegraphics[width=0.9\textwidth]{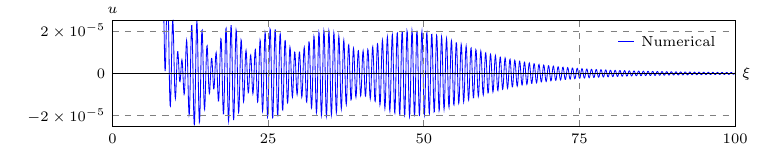}
 }
 
 \subfloat[Asymptotic contributions: pole (red) and saddle point (black)]{
\includegraphics[width=0.9\textwidth]{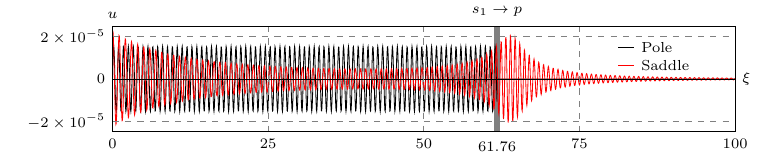}
 }
 
 \subfloat[Numerical solution (blue) and combined asymptotic contributions (black)]{
\includegraphics[width=0.9\textwidth]{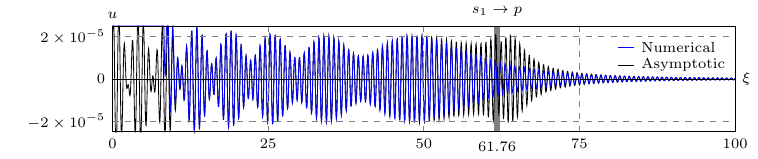}
 }
\caption{Comparison of numerical solution for $u(\xi,t)$ and asymptotic prediction of exponentially small pole and saddle contributions from the steepest descent method, for $c = 2$ and $\epsilon = 0.125$ at $t = 0.5$. The numerical solution is shown in blue in (a) and (b). In (c) the curves are the pole contribution (red) and the saddle point contribution (black). The point at which the pole contribution is switched off by the saddle point is shown as a grey line. In (d) the numerical solution (blue) is compared with the sum of the pole and saddle point contributions (black).}\label{fig:linnum}
\end{figure}

To validate the asymptotic saddle and pole contributions, we evalute the integral \eqref{e:kdvlinint} numerically along a contour that passes below the poles on the real axis.  A representative numerical example is shown in Figure \ref{fig:linnum} (a)--(b) for the case $c = 2$ with $\eps = 0.125$, evaluated at $t = 0.5$. Figure \ref{fig:linnum} (a) shows the shape of the full solution, while (b) depicts the exponentially small oscillations more clearly.

In Figure \ref{fig:linnum} (c), we have shown the asymptotic contributions from the pole and saddle point contributions, given by \eqref{eq:steadywaves} and \eqref{e:sdsaddle} respectively. The pole contributions are restricted to the region in Figure \ref{fig:ripple_region} that contains the stationary waves. In Figure \ref{fig:linnum} (d), we compare the combined asymptotic prediction with the numerical oscillations and see that the asymptotic analysis accurately describes the exponentially small oscillations.

The match is inaccurate in the neighbourhood of the Stokes switching point. This occurs as the saddle point approaches the pole location as $\xi \to 2\tau$ when $\tau$ becomes large, but the steepest descent analysis treats all of the contributions as separated in the complex plane. For completely descriptive asymptotics in this region, we would need to further rescale the problem to more accurately capture the behaviour in the neighbourhood of the saddle point/pole.

Also note that we treat the switching as instantaneous for the purposes of this analysis, but Stokes switching is a smooth process that takes place over a narrow region \cite{berry1988stokes, Daalhuis}; this does not have a significant impact on the asymptotic behaviour presented here.

\section{lKdV5: Factorial-Over-Power Analysis}\label{s:lKdV5}

The analysis from Section \eqref{s:linSD} was sufficient to completely determine the asymptotic behaviour of the small waves in the linearised KdV5 system governed by \eqref{e:lkdv5}; however, such an analysis is not possible for the nonlinear system governed by the full KdV5 equation \eqref{e:kdv5}. Nonlinear effects in space appear as convolutions in the transformed integral equation, which are typically more challenging to deal with using steepest descent analysis. Instead, we will use the method developed in \cite{Chapman}, which extends the exponential asymptotic ideas from \cite{Daalhuis} to study Stokes' phenomenon in nonlinear differential equations.

Throughout this section we use the variables $\sigma$, $\nu$, and $\alpha$ defined in \eqref{e:newvars}, and $p(\nu)$ defined in \eqref{e:p(nu)}.

\subsection{Series Expansion}

It is convenient to write \eqref{e:lkdv5} in terms of the analytically-continued shifted variable $\sigma = \xi + \i \pi/\sqrt{c}$ \eqref{e:newvars}, giving
\begin{align}\label{e:lkdv5sigma}
\eps^4 u_{\sigma\sigma\sigma\sigma\sigma} + \eps^2 u_{\sigma\sigma\sigma} + u_{\tau} = \eps^2 \left[  \tfrac{c}{2}\sech^2\left(\tfrac{\sqrt{c}}{2}\left(\sigma - \tfrac{\i\pi}{\sqrt{c}}\right)\right)\right]_{\sigma\sigma\sigma},& \qquad (\sigma,\tau) \in \mathbb{C} \times [0,\infty),
\end{align}
where 
\begin{equation}u(\sigma,0) = \tfrac{c}{2}\sech^2\left(\tfrac{\sqrt{c}}{2}\left(\sigma - \tfrac{\i\pi}{\sqrt{c}}\right)\right),
\end{equation}
and the far-field condition is still applied as $\xi \to \pm\infty$. We expand the solution as a power series in $\eps$, so
\begin{equation}\label{eq:useries}
u(\sigma,\tau) \sim \sum_{n=0}^{\infty}\eps^{2n} u^{(n)}(\sigma,\tau)\quad \as \quad \eps \to 0,
\end{equation}
where 
\begin{equation}\label{eq:uic}
u^{(0)}(\sigma,\tau) = \tfrac{c}{2}\sech^2\left(\tfrac{\sqrt{c}}{2}\left(\sigma - \tfrac{\i\pi}{\sqrt{c}}\right)\right).
\end{equation}
The recurrence relation obtained by balancing powers of $\eps$ at $\mathcal{O}(\eps^n)$ is given by
\begin{equation}\label{eq:lkdv5recur}
u^{(n-2)}_{\sigma\sigma\sigma\sigma\sigma} +  u^{(n-1)}_{\sigma\sigma\sigma}  + u^{(n)}_{\tau} = 0, \qquad u^{(n)}(\sigma,0) = 0,
\end{equation}
which is identical to \eqref{e:lkdv5recur} up to the shifted variable. 

\subsection{Late-Order Terms}\label{S:LOT_lkdv5}

To determine the effect of Stokes' phenomenon, we follow \cite{Chapman} and apply a late-order ansatz to the series terms to find the asymptotic behaviour of these terms in the limit $n \to \infty$, giving
\begin{equation}\label{e:ansatz}
u^{(n)} \sim \frac{U(\sigma,\tau)\Gamma(2n + \gamma)}{\chi(\sigma,\tau)^{2n+\gamma}}\quad \as \quad n \to \infty.
\end{equation}
By substituting this into \eqref{eq:lkdv5recur}, we obtain 
\begin{align}\nonumber
-\frac{\chi_{\sigma}^5 U\Gamma(2n + \gamma + 1)}{\chi^{2n+ \gamma+1}} +\frac{5\chi_{\sigma}^4 U_{\sigma}\Gamma(2n + \gamma)}{\chi^{2n+ \gamma}} + \frac{10\chi_{\sigma\sigma}\chi_{\sigma}^3 U \Gamma(2n + \gamma)}{\chi^{2n+ \gamma}}&\\ \nonumber
-\frac{\chi_{\sigma}^3 U\Gamma(2n + \gamma + 1)}{\chi^{2n+ \gamma+1}} +\frac{3\chi_{\sigma}^2 U_{\sigma}\Gamma(2n + \gamma)}{\chi^{2n+ \gamma}} + \frac{3\chi_{\sigma\sigma}\chi_{\sigma} U\Gamma(2n + \gamma)}{\chi^{2n+ \gamma}}&\\
-\frac{\chi_{\tau}U\Gamma(2n + \gamma + 1)}{\chi^{2n+ \gamma+1}} +\frac{U_{\tau}\Gamma(2n + \gamma)}{\chi^{2n+ \gamma}}& + \ldots = 0.
\end{align}
where the omitted terms are $o(u^{(n)})$ in the limit $n \to \infty$. The singulant equation is given by balancing at leading order in $n$ to obtain
\begin{equation}\label{e:singulant}
\chi_{\sigma}^5 + \chi_{\sigma}^3 + \chi_{\tau} = 0,
\end{equation}
with $\chi = 0$ at the location of the singularity under consideration. At the next order we obtain an equation for the prefactor $U$,
\begin{equation}\label{e:prefactor}
(5\chi_{\sigma}^4 + 3 \chi_{\sigma}^2) U_{\sigma} + U_{\tau} = -10 \chi_{\sigma}^3 \chi_{\sigma\sigma}U - 3 \chi_{\sigma}\chi_{\sigma\sigma} U.
\end{equation}
Determining the correct boundary condition for this equation typically requires an inner analysis in an appropriate neighbourhood of the singularity.

\subsubsection{Stationary terms}\label{s:t-independent}

We first identify the form of the late-order terms that is independent of $\tau$. This analysis is essentially identical to that of \cite{Grimshaw,trinh2010exponential}; we therefore omit the details, and only state the relevant results. We see from \eqref{eq:uic} that $u^{(0)}(\sigma,\tau)$ is singular at $\sigma = 0$. Solving \eqref{e:singulant} with this condition gives
\begin{equation}
\chi = \i\sigma.
\end{equation}
The prefactor equation \eqref{e:prefactor} reduces to $U_{\sigma} = 0$, so that $U$ is constant in $\sigma$; it is straight-forward to see from an inner expansion about $\sigma = 0$ that $U$ also does not vary in $\tau$. Finally, we determine $\gamma$ by comparing the strength of the singularity at $\sigma = 0$ from the late-order ansatz \eqref{e:ansatz} with the strength of the singularity in the leading-order solution. As the leading-order solution contains a double pole at $\sigma = 0$, consistency with the late-order terms requires $\gamma = 2$.

\subsubsection{Identifying the non-stationary singularity}\label{s:smallertime}

We must determine the location of the singularity associated with unsteady behaviour. This behaviour is not captured by the late-order terms in section \ref{s:t-independent}, which corresponds to the steady waves found in the steady solution. However, the appearance of this singularity is not as apparent as in \cite{chapman2005exponential}, where it could be seen in the first correction term.

Motivated by the observations in \cite{chapman2005exponential} about the form we expect this singularity to take, we calculate the asymptotic behaviour of $u(\sigma,\tau)$ in the neighbourhood of the initial disturbance at $\sigma = 0$ and $\tau =0$. While we are already in the boundary layer regime, we need to consider a new inner regime near this singular point in order to identify the form of this initial disturbance. We set $T = \eps \tau$ and $\eps^{1/4}X = \sigma$ in order for the time and space derivatives to remain present. This gives
\begin{equation}
u_{XXXXX} + \eps^{1/2}u_{XXX} + u_T = \eps^{1/2}u^{(0)}_{XXX}.
\end{equation}
We apply a WKB ansatz such that 
\begin{equation}\label{e:WKB1}
u(X,T) = A(X,T) \e^{\psi(X,T)/\eps},
\end{equation}
where $\phi = 0$ at $X = 0$. Applying \eqref{e:WKB1} and balancing the leading-order terms as $\eps \to 0$ gives
\begin{equation}
\psi_X^5 + \psi_T = 0.
\end{equation}
Solving the new differential equation gives the relevant solution in this regime as
\begin{equation}\label{e:smalltimelocal}
\psi = -\frac{4X^{5/4}}{5^{5/4}T^{1/4}}.
\end{equation}
This confirms that the solution $u$ is singular for small time, and that since both $X$ and $T$ appear with non-integer powers here, the relevant singular behaviour corresponds to $X = 0$ and $T = 0$. 

\subsubsection{Non-stationary terms}

The analysis of Section \ref{s:smallertime} suggests that we should solve \eqref{e:singulant} with the condition that $\chi = 0$ at $(\sigma,\tau) = (0,0)$. 
Solving this with Charpit's method requires writing the system of differential equations
\begin{align}
\diff{\sigma}{s} = 5 k_1^4 + 3 k_1^2,\qquad
\diff{\tau}{s} = 1,\qquad
\diff{\chi}{s} = 4k_1^5 + 2 k_1^3.
\end{align}
with $\sigma(0) = 0$, $t(0) = 0$ and $\chi(0) = 0$, in addition to the equation $k_1^5 + k_1^3 + k_2 = 0$. 
Integrating the differential equations and eliminating $k_1$, $k_2$, and $s$ gives
\begin{equation}
\chi = -\frac{\i\alpha^{1/2}(20\nu + \alpha)\tau}{25\sqrt{10}} =  p(\nu)\tau,
\end{equation} 
where $p$ is also given in \eqref{e:p(nu)}. 
The prefactor equation \eqref{e:prefactor} simplifies to
\begin{equation}
\frac{\sigma}{\tau}U_{\sigma} + U_{\tau} = -\frac{U}{2\tau}.
\end{equation}
which has solution
\begin{equation}\label{e:prefactorl}
U =\frac{B(\nu)}{\sqrt{\tau}},
\end{equation}
where $B$ is an unknown function of $\nu = \sigma/\tau$. Typically, determining $B(\nu)$ requires matching the late-order term ansatz with an inner expansion near the singularity as in \cite{chapman2005exponential}. This is a challenging computation even for linear problems, and is likely impossible for nonlinear problems (see Appendix \ref{a:inners}). The analysis in \cite{chapman2005exponential} uses the form of $B(\nu)$ explicitly in the second-generation Stokes line calculations; these Stokes lines appear due to new singularities in $B(\nu)$, and we can see from the steepest descents analysis in Section \ref{s:linSD} that such a singularity appears at $\nu = 2$. Nevertheless, we will leave $B(\nu)$ undetermined for now, and later utilise the condition that the initial condition is waveless in order to recover the necessary information. 

The late-order ansatz is therefore given by
\begin{equation}\label{e:ansatz2}
u^{(n)} \sim \frac{B(\nu)\Gamma(2n + \gamma)}{\sqrt{\tau}\chi^{2n+\gamma}}\quad \as \quad n \to \infty.
\end{equation}
To determine $\gamma$, we consider the behaviour of \eqref{e:ansatz2} in the limit $\tau \to 0$ and then $\sigma \to 0$, finding that
\begin{equation}
\chi \sim \frac{4 \i \sigma^{5/4}}{5^{5/4}\tau^{1/4}}\quad \as \quad \tau \to 0.
\end{equation}
If we assume that $B(\nu) \sim \nu^k$ in the limit that $\tau \to 0$ and then $\sigma \to 0$, we can enforce consistency of the late-order terms \eqref{e:ansatz2} with the leading order behaviour. In this sequence of limits, $u^{(0)}(\sigma,\tau) \sim -2/\sigma^2$. Balancing powers of $\tau$ and $\sigma$ between the late-order ansatz with $n$ taken to zero and this leading-order solution gives
\begin{equation}
k - \tfrac{5\gamma}{4} = -2, \qquad -k - \tfrac{1}{2} + \tfrac{\gamma}{4} = 0.
\end{equation} 
This tells us that $k = -1/8$ and $\gamma = 3/2$. We have therefore determined the late-order ansatz \eqref{e:ansatz} entirely, with the exception of the function $B(\nu)$.

\subsection{Stokes Lines}

Calculating $\chi$ associated with the stationary and non-stationary asymptotic contributions allows us to determine the location of Stokes lines in the solution, and to therefore reproduce the entire Stokes structure depicted in Figure \ref{fig:v_stokes}. The singulant associated with the stationary contribution is $\chi = \i\sigma$, while the singulant associated with the non-stationary contribution is $\chi = p(\nu)\tau$.

Stokes lines across which late-order terms associated with some exponential contribution $\chi_1$ switch on an exponential contribution associated with some $\chi_2$ satisfy
\begin{equation}
\mathrm{Im}[\chi_1] = \mathrm{Im}[\chi_2], \qquad \mathrm{Re}[\chi_1] < \mathrm{Re}[\chi_2],
\end{equation}
where the base power series is associated with $\chi = 0$. Comparing this condition with the equivalent condition in \eqref{e:SLsd} reveals the connection with the steepest descent analysis in Section \ref{s:linSD}. Similarly, the higher-order Stokes lines associated with three exponential contributions having $\chi_1$, $\chi_2$ and $\chi_3$ satisfy
\begin{equation}
\mathrm{Im}\left[\frac{\chi_2 - \chi_1}{\chi_3 - \chi_1}\right] = 0,
\end{equation}
which is equivalent to the condition given in \eqref{e:HOSLsd} for the steepest descent analysis. Hence, we have obtained enough information to completely reconstruct the Stokes lines and higher-order Stokes lines present in the asymptotic solution, which are depicted in Figure \ref{fig:v_stokes}. 

\subsection{Stokes Switching}\label{S:Switching_lkdv}

We have obtained the asymptotic form of the late-order terms \eqref{e:ansatz}, and can now use matched asymptotic expansions to determine the exponentially small quantity that appears as the Stokes lines are crossed \cite{Chapman,Daalhuis}. 

\subsubsection{Stationary contribution}

For the stationary late-order terms, the analysis is identical to \cite{trinh2010exponential}, and reveals that the quantity
\begin{equation}\label{e:steadywaves}
u_{\mathrm{exp}} \sim -\frac{2\pi\i}{\eps^2}\e^{\i\sigma/\eps} + \mathrm{c. c.} = -\frac{4\pi}{\eps^2}\e^{-\pi/(\sqrt{c}\eps)}\sin\left(\frac{\xi}{\eps}\right) \quad \mathrm{as} \quad \eps \to 0,
\end{equation}
is switched on as the associated Stokes line at $\mathrm{Re}(\xi) = 0$ is crossed from $\mathrm{Re}(\xi) < 0$ to $\mathrm{Re}(\xi) > 0$. This captures the appearance of the exponentially small ripples in the asymptotic solution, but not the time-dependent manner in which they propagate.

\subsubsection{Non-stationary contribution}\label{s:stokescalcs}

Motivated by the appearance of a characteristic variable $\nu$ in the prefactor, we rewrite the PDE \eqref{e:lkdv5sigma} giving
\begin{align}\label{e:lkdv5nu}
\eps^4 u_{\nu\nu\nu\nu\nu} + \eps^2 \tau^2 u_{\nu\nu\nu} -\nu \tau^4 u_{\nu} + \tau^5 u_{\tau} = \eps^2 \tau^2 \left[  \tfrac{1}{2}\sech^2\left(\tfrac{1}{2}\left(\nu\tau - \tfrac{\i\pi}{\sqrt{c}}\right)\right)\right]_{\nu\nu\nu},& \qquad (\nu,\tau) \in \mathbb{C} \times [0,\infty).
\end{align}
We then truncate the power series \eqref{eq:useries} after $N$ terms, obtaining
\begin{equation}\label{eq:useriesR}
u(\nu,\tau) = \sum_{n=0}^{N-1}\eps^{2n} u^{(n)}(\nu,\tau) + R(\nu,\tau),
\end{equation}
where $R$ is the exact remainder after truncation. Following the heuristic from \cite{Boyd2005}, we find that the optimal truncation is given by $N \sim |\chi|/(2\eps)$, and we therefore set $N = |\chi|/(2\eps) + \omega$, where $\omega \in [0,1)$ is chosen so that $N$ is an integer. Substituting \eqref{eq:useriesR} into \eqref{e:lkdv5nu} and simplifying the expression using \eqref{eq:lkdv5recur} gives
\begin{align}\label{eq:RNeqn}
\eps^4 R_{\nu\nu\nu\nu\nu} + \eps^2 \tau^2 R_{\nu\nu\nu} - \nu\tau^4 R_{\nu} + \tau^5 R_{\tau} =-\eps^{2N}\left[u^{(N-2)}_{\nu\nu\nu\nu\nu} + \tau^2 u^{(N-1)}_{\nu\nu\nu}\right]-\eps^{2N+2}u^{(N-1)}_{\nu\nu\nu\nu\nu}.
\end{align}
A WKB analysis of the homogeneous form of \eqref{eq:RNeqn} gives the solution $R \sim k U \e^{-\chi/\eps}$ as $\eps \to 0$ if the right-hand side is negligible, which is the case away from the Stokes lines. To capture the Stokes switching, we set $R = \mathcal{S}U\e^{-\chi/\eps}$. After using the form \eqref{e:ansatz} for the late-order terms on the right-hand side, \eqref{eq:RNeqn} becomes
\begin{equation}
\tau^5\pdiff{\mathcal{S}}{\tau}U\e^{-\chi/\eps} \sim \eps^{2N}\frac{\left(\chi_{\nu}^5 + \tau^2 \chi_{\nu}^3\right)U\Gamma(2N + \gamma +1)}{\chi^{2N + \gamma + 1}}+\eps^{2N+2}\frac{\chi_{\nu}^5 U\Gamma(2N + \gamma +3)}{\chi^{2N + \gamma + 3}},
\end{equation}
where the omitted terms are smaller in the limit $\eps \to 0$ than those which remain. Recalling that $\chi = p(\nu)\tau$ and writing the partial derivative in terms of $\chi$ gives, after some simplification,
\begin{equation}
\pdiff{\mathcal{S}}{\chi} \sim \left[\eps^{2N}\frac{\left(p_{\nu}^5 + p_{\nu}^3\right)\Gamma(2N + \gamma +1)}{\chi^{2N + \gamma + 1}}+\eps^{2N+2}\frac{p_{\nu}^5 \Gamma(2N + \gamma +3)}{\chi^{2N + \gamma + 3}}\right]\frac{\e^{\chi/\eps}}{p}.
\end{equation}
Motivated by the form of $N$, we let $\chi = r\e^{\i\theta}$ and consider the variation in the $\theta$-direction, which is the direction in which the change is rapid. After some simplification, we find
\begin{equation}
\pdiff{\mathcal{S}}{\theta} \sim \frac{\i}{p \eps^{\gamma}} \sqrt{\frac{2\pi r}{\eps}}\left[(p_{\nu}^5 + p_{\nu}^3)\e^{-\i\theta} + p_{\nu}^5 \e^{-3\i\theta}\right]\exp\left(\frac{r\e^{\i\theta}}{\eps} - \frac{r}{\eps} - \frac{\i\theta r}{\eps} - \i \theta(\omega + \gamma)\right)
\end{equation}
Letting $\vartheta = \eps^{1/2}\theta$ and integrating the expression gives the jump across the Stokes line as
\begin{align}
\left[\mathcal{S}\right]_-^+ &\sim \left(\frac{2p_{\nu}^5 + p_{\nu}^3}{p}\right)\frac{2\pi\i}{\eps^{\gamma}}=\frac{\pi\i}{\eps^{\gamma}},
\end{align}
We have therefore determined the asymptotic form of the time-dependent exponentially small contribution,
\begin{equation}\label{e:u_exp}
R \sim \frac{\pi\i B(\nu)}{\tau^{1/2}\eps^{3/2}}\e^{-\chi/\eps} \quad \mathrm{as} \quad \eps \to 0,
\end{equation}
which is active in the region $\Im[p]> 0$. 

While it is not necessary for our analysis, we can now identify $B(\nu)$ by matching the exponentially small term \eqref{e:u_exp} with the corresponding steepest descent term $I_{s1}$ from \eqref{e:sdsaddle}. This gives
\begin{equation}
B(\nu) =\frac{-\i 2^{3/4}\alpha^{5/4}\e^{\i\pi/4}}{5^{1/4}\sqrt{\pi}(\alpha - 10)\sqrt{\alpha - 3}}.
\end{equation}

\subsection{Second Generation Late-Order Terms}

The analysis in Section \ref{S:Switching_lkdv} showed where the two relevant exponential contributions in the solution appear, but this cannot be the entire Stokes structure of the asymptotic solution, as it does not satisfy condition that there are no waves at $\tau = 0$, nor does it satisfy the far-field conditions; the exponentially small oscillations \eqref{e:steadywaves} are present to the right of $\mathrm{Re}[\sigma] = 0$ and continue as $\xi \to \infty$ with constant, nonzero amplitude. 

In order for the far-field conditions to be satisfied, the solution must contain an additional Stokes line across which the time-varying exponential contribution causes the oscillatory waves to disappear. This Stokes line is caused by the divergence of the asymptotic series for $R$, which causes the appearance of new exponentially small terms and is sometimes known as second-generation Stokes switching. 

To determine the exponentially small correction to $R$, we write it as an exponentially scaled power series as $\eps \to 0$, 
\begin{equation}\label{e:Rseries}
R \sim \e^{-p\tau/\eps}\sum_{n=0}^{\infty}A^{(n)}(\nu,\tau)\eps^{n-3/2}\quad \as \quad \eps \to 0.
\end{equation}
From \eqref{e:u_exp}, we know that 
\begin{equation}\label{e:2genf0}
A^{(0)} = \frac{2^{3/4}\alpha^{5/4}\e^{\pi\i/4}\sqrt{\pi}}{5^{1/4}(\alpha-3)^{1/2}(10-\alpha)\tau^{1/2}} \sim -\frac{2\e^{\pi \i/4}\sqrt{7\pi}}{\tau^{1/2}(\nu-2)}\quad \as \quad \nu \to 2.
\end{equation}
The singularity in $A^{(0)}$ at $\nu = 2$ causes divergence that produces its own new Stokes lines. This is apparent from the algebraic expression in \eqref{e:2genf0}; however, we can also see this because we expect $\chi_1 = p(\nu)\tau$ to switch on $\chi_2 = \i\nu\tau$. It is simple to show that the equation $\chi_2 - \chi_1 = 0$, or
\begin{equation}\label{e:magic_equation_1}
\i \nu - p = 0
\end{equation}
is solved by $\nu = 2$. Hence, if we are able to assume that the waves do not radiate infinitely quickly, we do not need to use the exact form for $B(\nu)$. We will instead assume that
\begin{equation}\label{e:kappa}
A^{(0)}(\nu) \sim \frac{C}{\tau^{1/2}}(\nu - 2)^{\kappa} \quad \mathrm{as} \quad \nu \to 2,
\end{equation}
where the scaling of $\tau$ is due to the leading-order exponentially small behaviour from \eqref{e:u_exp}, and $C$ and $\kappa$ are not yet determined. We require $\kappa$ to not be a positive integer so that $A^{(0)}(\nu)$ is singular at $\nu = 2$.  We note from \eqref{e:2genf0} that $A^{(0)}$ is also singular at $\nu = -9/20$, but any contributions that arise due to this are not relevant to the current switching analysis. 

We substitute \eqref{e:Rseries} into the homogeneous equation for $R$, as we expect this switching to occur away from the regular Stokes line. This equation is
\begin{equation}\label{e:2genReq}
\eps^4 R_{\nu\nu\nu\nu\nu} + \eps^2 \tau^2 R_{\nu\nu\nu} - \nu \tau^4 R_{\nu} + \tau^5 R_{\tau} = 0.
\end{equation}
Substituting in \eqref{e:Rseries} and balancing terms at $\mathcal{O}(\e^{-p\nu/\eps}\eps^{n-3/2})$ gives a recurrence relation for $A^{(n)}$:
\begin{align}
\Big(5 p_{\nu} \tau A_{\nu}^{(n)} - &10 p_{\nu}^3 \tau^3 A_{\nu\nu}^{(n-1)} + 10 p_{\nu}^2 \tau^2 A_{\nu\nu\nu}^{(n-2)} - 5 p_{\nu} \tau A_{\nu\nu\nu\nu}^{(n-3)} + A_{\nu\nu\nu\nu\nu}^{(n-4)} + 10 p_{\nu}^3 p_{\nu\nu}  \tau^4 A^{(n)} 
\nonumber\\& 
- 30 p_{\nu}^2 p_{\nu\nu}^2 \tau^4 A_{\nu}^{(n-1)}
+30 p_{\nu}  p_{\nu\nu} \tau^2 A_{\nu\nu}^{(n-2)} - 10 p_{\nu\nu} \tau A_{\nu\nu\nu}^{(n-3)}\Big)+ \tau^2\Big(3p_{\nu} \tau A_{\nu}^{(n)} 
\nonumber\\ &
- 3  p_{\nu} \tau A_{\nu\nu}^{(n-1)} + A_{\nu\nu\nu}^{(n-2)} - 3 p_{\nu} p_{\nu\nu} \tau^2 A^{(n)} - 3 p_{\nu\nu} \tau A_{\nu}^{(n-1)}\Big) - \nu \tau^4 A_{\nu}^{(n)} + \tau^5 A_{\tau}^{(n)} = 0.\label{e:2genArecur}
\end{align}
As $\chi = p(\nu)\tau$, a balancing argument to determine the dependence of $A^{(n)}$ on $\tau$ gives
\begin{equation}\label{e:f0local}
A^{(n)} = \frac{f^{(n)}(\nu)}{\tau^{n + 1/2}}, \qquad f^{(0)}(\nu) \sim C(\nu-2)^{\kappa} \quad \mathrm{as} \quad \nu \to 2,
\end{equation}
where the offset power of $\tau$ in $A^{(n)}$ is chosen to be $-1/2$ so that $f^{(0)}(\nu)$ is consistent with $A^{(0)}$ from \eqref{e:kappa}. 
Substituting this into \eqref{e:2genArecur} and using $\chi = p(\nu)\tau$ gives the ordinary differential-difference equation
\begin{align}
\Big(5 p_{\nu}^4 f^{(n)}_{\nu} - 10 p_{\nu}^3 f^{({n-1})}_{\nu\nu} + 10 p_{\nu}^2 f^{({n-2})}_{\nu\nu\nu} - 5 p_{\nu} f^{({n-3})}_{\nu\nu\nu\nu}
+ f^{({n-4})}_{\nu\nu\nu\nu\nu}\Big)&\nonumber \\+ \Big(3 p_{\nu}^2  f^{(n)}_{\nu} - 3 p_{\nu} f^{({n-1})}_{\nu\nu} + f^{({n-2})}_{\nu\nu}\Big) - \nu f^{(n)}_{\nu}& - (n + 1/2) f^{(n)} + \ldots = 0,\label{e:frecur}
\end{align}
where the omitted terms are subdominant in the limit $n \to \infty$ compared to those retained. We can now use late-order term analysis to determine the behaviour of $f^{(n)}$ for large $n$.
As in the previous section, we assume a factorial-over-power late-order ansatz for $f^{(n)}$, given by
\begin{equation}\label{e:lateF}
f^{(n)} \sim \frac{F(\nu)\Gamma(n + \beta)}{\mu(\nu)^{n+\beta}}\quad \as \quad n \to \infty,
\end{equation}
with $\mu = 0$ at $\nu = 2$. This implies that the late-order terms for in the series expansion for $R$ from \eqref{e:Rseries} have the factorial-over-power form
\begin{equation}\label{e:lateA}
A^{(n)} \sim \frac{F(\nu)\Gamma(n + \beta)}{\phi(\nu,\tau)^{n+\beta}}\quad \as \quad n \to \infty,
\end{equation}
where $\phi(\nu,\tau) = \mu(\nu)\tau$. 

Substituting the late-order ansatz \eqref{e:lateF} into the recurrence relation \eqref{e:frecur} allows us to determine the prefactor and singulant in a similar fashion to the analysis performed in Section \ref{S:LOT_lkdv5}. Taking the leading-order solution in the limit $n \to \infty$ gives the singulant equation
\begin{equation}
\mu + \mu_{\nu}^5 + 5 p_{\nu}\mu_{\nu}^4 + (10 p_{\nu}^2 + 1)\mu_{\nu}^3 + (10 p_{\nu}^2 + 3 )p_{\nu}\mu_{\nu}^2 = 0,
\end{equation}
with the condition $\mu = 0$ at $\nu = 2$. This could be solved using Charpit's method, but direct substitution is sufficient to confirm that the solution is 
\begin{equation}\label{e:fsingulant}
\mu = \i \nu - p,
\end{equation} 
as expected. Balancing the result of substituting \eqref{e:lateF} into \eqref{e:frecur} to the next order in the limit $n \to \infty$ gives an equation for the prefactor $F(\nu)$.
Using \eqref{e:fsingulant} to obtain $\mu_{\nu} = \i - p_{\nu}$ and $\mu_{\nu\nu} = -p_{\nu\nu}$ allows us to simplify this expression, resulting in
\begin{equation}
(\beta - 1/2)F + (2-\nu)F_{\nu} = 0.
\end{equation}
This can be solved directly to give
\begin{equation}\label{e:Feq0}
F = \Lambda(\nu-2)^{\beta-1/2},
\end{equation}
where $\Lambda$ is an arbitrary constant. In order to determine $\beta$, we must consider the behavior of the late-order ansatz \eqref{e:lateF} in the neighbourhood of the singularity at $\nu = 2$, and require it to be consistent with the behaviour of $f^{(0)}$ in this region.
In the limit $\nu \to 2$, we have
\begin{equation}\label{e:2genlocal}
\mu \sim -\frac{\i}{28}(\nu-2)^2,
\end{equation}
This means that the local behaviour of the late-order ansatz \eqref{e:lateF} is given by
\begin{equation}\label{e:lateFlocal}
f^{(n)} \sim \frac{\Lambda(\nu-2)^{\beta-1/2}\Gamma(n+\beta)}{(-\i (\nu-2)^2/28)^{n+\beta}}\quad \mathrm{as} \quad n \to \infty,\, \nu \to 2.
\end{equation}
For \eqref{e:lateF} to be compatible with $f^{(0)}$, given in \eqref{e:f0local}, in this limit, the local late-order ansatz \eqref{e:lateFlocal} must have the same singularity strength at $\nu = 2$ as $f^{(0)}$ when $n$ is set to 0 in the late-order ansatz. This requires that $\beta - 1/2 - 2\beta = \kappa$, or $\beta = -1/2 - \kappa$. Hence
\begin{equation}\label{e:fpref}
F = \Lambda (\nu-2)^{\kappa}.
\end{equation}

Finally, we can determine $\Lambda$ and $\kappa$ by matching the late-order terms with the solution in an inner region near $\nu = 2$. To obtain the appropriate scaling, we first set 
\begin{equation}\label{e:scaling}
R = A(\nu,\tau)\eps^{-3/2}\e^{-p\tau/\eps}.
\end{equation}
Substituting this directly into the equation for $R$ from \eqref{e:2genReq} gives 
\begin{align}
\Big(3\tau^4 p_{\nu}p_{\nu\nu} + 10 \tau^4 p_{\nu}^3 p_{\nu\nu} A + \tau^5A_{\tau} - \tau^4 \nu A_{\nu} + 3 \tau^4 p_{\nu}^2 A_{\nu} + 5 \tau^4 p_{\nu}^4 A_{\nu}\Big)
- \eps\Big(15\tau^3 p_{\nu}p_{\nu\nu}A + \tau^3 p_{\nu\nu\nu}A 
\nonumber&\\ 
+10 \tau^3 p_{\nu}^2 p_{\nu\nu\nu}A + 3 \tau^3 p_{\nu\nu}A_{\nu} 
+30 \tau^3p_{\nu}^2 p_{\nu\nu}A_{\nu} + 3 \tau^3 p_{\nu}A_{\nu\nu} + 10\tau^3 p_{\nu}^3 A_{\nu\nu}\Big)
+\eps^2\Big(10\tau^2 p_{\nu\nu}p_{\nu\nu\nu}A 
\nonumber&\\ 
+ 5 \tau^2 p_{\nu}p_{\nu\nu\nu\nu}A 
+15 \tau^2 p_{\nu\nu}^2 A_{\nu} + 20 \tau^2 p_{\nu}p_{\nu\nu\nu}A_{\nu} 
+ 30 \tau^2 p_{\nu}p_{\nu\nu}A_{\nu\nu} + \tau^2 A_{\nu\nu\nu} + 10 \tau^2 p_{\nu}^2 A_{\nu\nu\nu}\Big)
\nonumber&\\
-\eps^3\Big(\tau p_{\nu\nu\nu\nu\nu}A+5 \tau p_{\nu\nu\nu\nu}A_{\nu}+ 10 \tau p_{\nu\nu\nu}A_{\nu\nu} + 10 \tau p_{\nu\nu} A_{\nu\nu\nu} - 5 \tau p_{\nu}A_{\nu\nu\nu\nu}\Big) + \eps^4 A_{\nu\nu\nu\nu\nu} &= 0.\label{e:2geninnerA}
\end{align}
From \eqref{e:2genlocal} and recalling from \eqref{e:lateA} that $\phi(\nu,\tau) = \mu(\nu)\tau$, we see that $\phi \sim -\tfrac{\i}{28}(\nu-2)^2 \tau$ near $\nu = 2$. From the expression in \eqref{e:lateFlocal}, we determine that the series for $R$ in \eqref{e:Rseries} fails to be asymptotic if $\phi = \mathcal{O}(\eps)$. Hence, we define a new inner coordinate $y$ such that
\begin{equation}\label{e:2geny}
\nu  = 2 + \eps^{1/2}\tau^{-1/2}y.
\end{equation}
Writing \eqref{e:2geninnerA} in terms of this new variable and taking the leading order in the limit $\eps \to 0$ gives
\begin{equation}\label{e:2geninner}
\frac{\tau^4}{2}(A - 2\tau A_{\tau} + y A_y + 14 \i A_{yy}) = 0.
\end{equation}
We define a series expansion for $A$ in this inner region,
\begin{equation}\label{e:2geninnerAser}
A \sim \sum_{n=0}^{\infty} \frac{a^{(n)}(\tau)}{y^{2n+1}} \quad \as \quad |y|\to \infty.
\end{equation}
Substituting \eqref{e:2geninnerAser} into \eqref{e:2geninner} gives
\begin{equation}
a^{(n)} + 2\tau  a^{(n)}_{\tau} - (2n+1)a^{(n)} + 14 \i(2n-1)(2n)a^{(n-1)} = 0,
\end{equation}
This simplifies to $a^{(n)} = 28\i(n-1/2)a^{(n-1)}$, which is solved by 
\begin{equation}\label{eq:innergamma}
a^{(n)} = \frac{a^{(0)} (28\i)^n \Gamma(n+1/2)}{\sqrt{\pi}}.
\end{equation}
Writing the late-order terms from \eqref{e:lateA} using the local behaviour \eqref{e:lateFlocal} in inner coordinates \eqref{e:2geny} gives
\begin{equation}
\eps^{n+1/2}A^{(n)} \sim \frac{\Lambda\Gamma(n+1/2)(28\i)^{n+1/2}}{y^{2n+1}} \quad \mathrm{as} \quad n \to \infty,\, \nu \to 2.
\end{equation}
Balancing this with the corresponding term from the inner series gives
\begin{equation}\label{e:lambdabeta}
\Lambda = \frac{a^{(0)}}{\sqrt{\pi}(28\i)^{1/2}}, \qquad \beta = \frac{1}{2}.
\end{equation}
Hence, we also find that $\kappa = 0$ and $F = \Lambda$, which is consistent with the result from the steepest descent analysis \eqref{e:2genf0}. The one quantity that remains to be determined is the constant $a^{(0)}$.

We note here that $a^{(0)}$ can be calculated precisely using the local behaviour of $f^{(0)}$ near $\nu = 2$, which was found in the steepest descent analysis. The value of $A^{(0)}$ in \eqref{e:2genf0} can be used to calculate $f^{(0)}$ and subsequently determine that $a^{(0)} = -2\e^{\i\pi/4}\sqrt{7\pi}$. The result in \eqref{e:lambdabeta} then gives $\Lambda = -1$. We will later obtain this result without reference to the steepest descent analysis.

\subsection{Second Generation Stokes Switching}\label{s:2gen_lkdv5}
The analysis in this section follows largely identical steps to Section \ref{S:Switching_lkdv}, with some additional algebraic terms that do not significantly change the process. We will therefore omit the algebraic details. 

We first truncate the series for $R$ from \eqref{e:Rseries} after $M$ terms, where we will choose $M$ such that the truncation is optimal. This gives
\begin{equation}\label{e:2gentruncated}
R = \frac{1}{\eps^{3/2}}\left(\sum_{n=0}^{M-1}\eps^n A^{(n)} + Q\right) \e^{-p\tau/\eps},
\end{equation}
where $Q$ is the truncation remainder. We then substitute this equation into the governing equation for $R$ from \eqref{e:2genReq}. Using \eqref{e:2genArecur} to cancel terms, it is possible to obtain a direct equation for $Q$. 
By setting $Q = \mathcal{S}\Lambda\e^{-\mu \tau/\eps}$, this simplifies to
\begin{align}
\tau^5 \pdiff{\mathcal{S}}{\tau} \Lambda \e^{-\mu \tau/\eps} \sim \, & \eps^{M+3}A_{\nu\nu\nu\nu\nu}^{(M-1)} + \eps^{M+2}\Big[5\tau p_{\nu}A_{\nu\nu\nu\nu}^{(M-1)} + A_{\nu\nu\nu\nu\nu}^{(M-2)}\Big] 
\nonumber \\ & 
+ \eps^{M+1}\Big[(10\tau^2 p_{\nu}^2  \tau^2)A_{\nu\nu\nu}^{(M-1)}  + 5 \tau p_{\nu} A_{\nu\nu\nu\nu}^{(M-2)} + A_{\nu\nu\nu\nu\nu}^{(M-3)}\Big]
\nonumber \\ & 
+ \eps^{M}\Big[(10 \tau^2 p_{\nu}^3 + 3 \tau^2 p_{\nu}) A_{\nu\nu}^{(M-1)} + (10\tau^2 p_{\nu}^2  \tau^2)A_{\nu\nu\nu}^{(M-2)}  + 5 \tau p_{\nu} A_{\nu\nu\nu\nu}^{(M-3)} + A_{\nu\nu\nu\nu\nu}^{(M-4)}\Big],
\end{align}
Applying a nearly identical sequence of steps to those outlined in Section \ref{S:Switching_lkdv} reveals that the quantity appearing across the second-generation Stokes line is
\begin{equation}
[\mathcal{S}]_-^+ \sim -\frac{2\pi\i}{\eps^{1/2}}\quad \mathrm{as} \quad \eps \to 0.
\end{equation}
Recalling that $\mu = \i\nu - p$ and $Q$ is the truncation remainder of $R$ from \eqref{e:2gentruncated}, this means that the asymptotic contribution $Q$, where
\begin{equation}\label{e:newoscillations}
Q \sim -\frac{2\pi\i\Lambda}{\eps^{2}} \e^{-\i\nu\tau/\eps} \quad \mathrm{as} \quad \eps \to 0
\end{equation}
is present in the expansion for $R$, and hence $u$, in the region $\mathrm{Im}(\i\nu - p) > 0$. As identified in the steepest descent analysis, this region includes both the initial data on $\tau = 0$ and the far-field $\xi \to \infty$, where the solution is known to be waveless. Hence, the exponentially small oscillations generated by this Stokes line must cancel precisely with the stationary late-order oscillations from \eqref{e:steadywaves}. This cancellation requires that $\Lambda = -1$ (and hence that $a^{(0)} = -2\e^{\i\pi/4}\sqrt{7\pi}$). 

By selecting this value of $\Lambda$, we are able to determine the exponentially small behaviour that satisfies the condition that $u(\xi,0)$ does not contain any exponentially small oscillations. We have been able to completely recover the wave behaviour found in Section \ref{s:linSD}, and hence the information presented in Figure \ref{fig:v_stokes} using the method of steepest descents.

\subsection{The need for $B(\nu)$}\label{s:Bnu}

As discussed previously, determining the explicit form of $B(\nu)$ through an inner analysis in the neighbourhood of the singularity is challenging, whether using the integral formulation of \cite{chapman2005exponential} or the direct late-order analysis of \cite{shelton2025exponential}. For nonlinear problems, such a calculation is likely to be intractable. While numerical values of $B(\nu)$ can be obtained for individual values of $\nu$, it is difficult to extract from these computations the information needed to complete the asymptotic analysis.

In the preceding sections, there were three instances in which the form of $B(\nu)$ would normally be required. The relevant information, and the way it was obtained without explicitly calculating $B(\nu)$, are summarised below.
\begin{itemize}
\item The location of the singularity of $B(\nu)$ at $\nu=2$, which corresponds to the emergence of a new Stokes line across which the stationary waves are switched. Rather than determining this singularity directly, we required that the exponential contribution associated with the new Stokes line have the same phase as the existing stationary waves so that the two asymptotic contributions can interact. Solving \eqref{e:magic_equation_1} under this condition yields $\nu=2$.

\item The value of $\beta$ in the late-order ansatz \eqref{e:lateF}, and hence the strength of the singularity of $B(\nu)$ at $\nu=2$, denoted by $\kappa$ in \eqref{e:kappa}. This value was obtained by matching the late-order ansatz \eqref{e:lateF} to the local expansion near $\nu=2$ given by \eqref{eq:innergamma}. The matching is only possible when $\beta=1/2$.

\item The value of $\Lambda$, and hence $a^{(0)}$. This quantity is normally determined from the coefficient $C$ in \eqref{e:kappa}, which multiplies the dominant singular behaviour of $B(\nu)$ near $\nu=2$. Instead, we matched the amplitude of the stationary oscillations \eqref{e:steadywaves} with the oscillations generated across the new Stokes line \eqref{e:newoscillations}. Requiring these contributions to cancel so that the solution is waveless as $u \to \infty$, and also at $\tau = 0$, gives the result that $\Lambda=-1$.
\end{itemize}

In each case, the required information was obtained by enforcing the condition that $u$ is waveless at $\tau = 0$, which also implies that $u\to0$ as $\xi\to\infty$. The only way this could occur is if the contribution appearing across the second-generation Stokes line precisely cancels the stationary waves, which in turn provides the necessary information about $B(\nu)$.

 \section{KdV5: Late-Order Term Analysis}\label{s:KdV5}
 
Finally, we return to the nonlinear KdV5 problem, \eqref{e:kdv5ST}--\eqref{e:kdv5ST-ff}. In this section, we will demonstrate that the exponential asymptotic analysis is essentially identical to that of Section \ref{s:lKdV5}, and that the only significant change is the value of $B(\nu)$ that we obtain. 

The stationary wave contribution for KdV5 was calculated in \cite{Grimshaw,trinh2010exponential}. These studies found that the exponentially small stationary oscillations are present for $\xi > 0$ and have the form
\begin{equation}\label{e:kdv5steadywaves}
u_{\mathrm{exp}} \sim -\frac{2\pi\i\Lambda}{\eps^2}\e^{\i\sigma/\eps} + \mathrm{c. c.} = -\frac{4\pi\Lambda}{\eps^2}\e^{-\pi/(\sqrt{c}\eps)}\sin\left(\frac{\xi}{\eps}\right) \quad \mathrm{as} \quad \eps \to 0,
\end{equation}
where $\Lambda \approx -19.97$. 

\subsection{Series Expansion}

To determine the contribution associated with the non-stationary exponentially small terms, we perform a late-order term asymptotic analysis. We expand the solution as a power series of the form \eqref{eq:useries}, which has the leading-order solution \eqref{eq:uic}. The recurrence relation for the power series terms is now given by
\begin{align}\label{e:kdv5recurLOT}
w^{(n-2)}_{\sigma\sigma\sigma\sigma\sigma} + w^{(n-1)}_{\sigma\sigma\sigma} + 6 \sum_{j=0}^{n-1}w^{(j)} w^{(n-1-j)}_{\sigma}-cw^{(n-1)}_{\sigma}+  w^{(n)}_{\tau} = 0,& \qquad (\xi,\tau) \in (-\infty,\infty) \times [0,\infty),
\end{align}
which is the same as \eqref{e:kdv5recur} up to the shifted variable. We assume a late-order ansatz of the form \eqref{e:ansatz} to describe $u^{(n)}$ in the limit $n \to \infty$, and substitute this into \eqref{e:kdv5recurLOT}. We now arrive at the rationale for constructing the auxiliary forced linear PDE \eqref{e:lkdv5}. The linear equation was designed so that once the ansatz \eqref{e:ansatz} was applied, the two leading orders as $n \to \infty$ contain exactly the same terms. 

\subsection{Late-Order Terms and Stokes Lines}

Balancing the leading-order terms gives the singulant and prefactor equations
\begin{equation}\label{e:singulantandprefacto}
\chi_{\sigma}^5 + \chi_{\sigma}^3 + \chi_{\tau} = 0, \qquad (5\chi_{\sigma}^4 + 3 \chi_{\sigma}^2) U_{\sigma} + U_{\tau} = -10 \chi_{\sigma}^3 \chi_{\sigma\sigma}U - 3 \chi_{\sigma}\chi_{\sigma\sigma} U,
\end{equation}
which are identical to \eqref{e:singulant} and \eqref{e:prefactor}, as expected. The analysis of the late-order terms remains the same and the singulants are therefore identical, so
\begin{equation}\label{e:singulantskdv5}
\chi = p(\nu)\tau,
\end{equation} 
where $\nu = \sigma/\tau$, as well as the corresponding complex conjugate solutions. We therefore expect all of the Stokes line interactions to remain the same as the lKdV5 analysis, and hence the Stokes structure is identical to that given in Figure \ref{fig:v_stokes}, even though a steepest descent analysis is not possible for this nonlinear problem. The corresponding prefactor calculations give
\begin{equation}\label{e:prefactorskdv5}
U = \frac{\tilde{B}(\nu)}{\sqrt{\tau}},
\end{equation}
and we accordingly find $\gamma = 3/2$. The quantity $\tilde{B}(\nu)$ is an arbitrary function of $\nu$; we use the tilde to emphasise that this is a different quantity to $B(\nu)$ from \eqref{e:prefactorl}. Aside from the behaviour of $\tilde{B}(\nu)$, these results are all identical to the lKdV5 analysis. 

Typically, $\tilde{B}(\nu)$ would be calculated as in \cite{chapman2005exponential} by scaling the problem to consider a local expansion near singularities in the analytically-continued solution. In this local expansion the PDE is rescaled such that the nonlinear effects contribute at all orders of the solution, which means that determining an exact form for the expansion is likely to be challenging or intractable. For further discussion of how this inner expansion could be approached, see Appendix \ref{a:inners}.

As discussed in Section \ref{s:Bnu}, however, we can avoid calculating this inner expansion entirely by instead enforcing the condition that the solution has no oscillations present at $\tau = 0$, and therefore the stationary contribution and the second-generation contribution must cancel precisely to permit this to occur.

Writing the PDE in terms of the variables $\nu$ and $\tau$ gives
\begin{equation}\label{e:KdV5nu}
\eps^4 w_{\nu\nu\nu\nu\nu} + \eps^2 \tau^2 w_{\nu\nu\nu} + 6 \tau^4 w w_{\nu} - c\tau^4 w_{\nu} - \nu \tau^4 w_{\nu} + \tau^4 w_{\tau} = 0.
\end{equation}
By truncating the power series optimally after $N$ terms, we obtain a remainder equation
\begin{align}\nonumber
\eps^4 R_{\nu\nu\nu\nu\nu}  + \eps^2 \tau^2 R_{\nu\nu\nu} + 6  \tau^4 R \sum_{j=0}^{N-1}\eps^{2j+2}w^{(j)}_{\nu} + 6  \tau^4 R_{\nu} \sum_{j=0}^{N-1}\eps^{2j+2}w^{(j)}+ 6  \eps^2\tau^4 R R_{\nu} - \eps^2  c \tau^4 R_{\nu} - \nu \tau^4 R_{\nu}& \\
+ \,\tau^4 w_{\tau}  = -\eps^{2N}\left[w^{(N-2)}_{\nu\nu\nu\nu\nu} + \tau^2 w^{(N-1)}_{\nu\nu\nu}\right]-\eps^{2N+2}w^{(N-1)}_{\nu\nu\nu\nu\nu} - 6\tau^4 \sum_{j=N}^{2N-2}\eps^{2j+2}\sum_{\ell = j - (N-1)}^{N-1} w^{(j-\ell)}w_{\nu}^{(\ell)}.&
\label{e:KdV5R}
\end{align}
Retaining only the leading order behavior as $\eps \to 0$ (and hence $N \to \infty$) gives the same remainder equation as in \eqref{eq:RNeqn}. The subsequent analysis proceeds identically to Section \ref{s:stokescalcs}, meaning that the associated exponential contribution is identical to \eqref{e:u_exp}, and hence given by 
\begin{equation}\label{e:kdv5u_exp}
R \sim \frac{\pi\i \tilde{B}(\nu)}{\tau^{1/2}\eps^{3/2}}\e^{-\chi/\eps} \quad \mathrm{as} \quad \eps \to 0,
\end{equation}
which is active in the region $\Im[p]> 0$. The only difference between this expression and \eqref{e:u_exp} is the form of $\tilde{B}(\nu)$. We were able to calculate the equivalent term, $B(\nu)$, for the lKdV5 analysis using steepest descents, but this option is not available for the nonlinear KdV5 analysis. We expect that the nonlinear contributions would make such a calculation intractable, and therefore cannot expect to make use of a convenient analytical form for $\tilde{B}(\nu)$. 

\subsection{Second-Generation Stokes Switching}

As in the lKdV5 analysis, the next step is to determine a series expansion for $R$, where \eqref{e:kdv5u_exp} is the leading-order behaviour. We therefore apply the expansion \eqref{e:Rseries}, noting that we do not have an expression like \eqref{e:2genf0} for $A^{(0)}$.  We are able to conclude that $\tilde{B}(\nu)$ must have a singularity at $\nu = 2$ by the same reasoning that produced \eqref{e:magic_equation_1}, and hence may again assume the form \eqref{e:kappa}. The homogeneous remainder equation is now given by
\begin{align}\nonumber
\eps^4 R_{\nu\nu\nu\nu\nu}  + \eps^2 \tau^2 R_{\nu\nu\nu} + 6  \tau^4 R \sum_{j=0}^{N-1}\eps^{2j+2}w^{(j)}_{\nu} + 6 \tau^4 &R_{\nu} \sum_{j=0}^{N-1}\eps^{2j+2}w^{(j)}\\+ 6  \eps^2\tau^4 R R_{\nu}& -\eps^2 c \tau^4   R_{\nu} - \nu \tau^4 R_{\nu} + \tau^5 R_{\tau}  = 0
\label{e:KdV5R_hom}
\end{align}
Substituting in \eqref{e:Rseries} and balancing terms at $\mathcal{O}(\e^{-p \tau/\eps} \eps^{n-3/2})$ gives a recurrence relation for $A^{(n)}$ that is similar to \eqref{e:2genArecur_kdv5}, but with the addition of further terms. The recurrence relation for $A^{(n)}$ is given by
\begin{align}\nonumber
\Big(5 p_{\nu} \tau A_{\nu}^{(n)} - 10 p_{\nu}^3 \tau^3 A_{\nu\nu}^{(n-1)} + 10 p_{\nu}^2 \tau^2 A_{\nu\nu\nu}^{(n-2)} - 5 p_{\nu} \tau A_{\nu\nu\nu\nu}^{(n-3)} + A_{\nu\nu\nu\nu\nu}^{(n-4)}& \\
+ 10 p_{\nu}^3 p_{\nu\nu}  \tau^4 A^{(n)} - 30 p_{\nu}^2 p_{\nu\nu}^2 \tau^4 A_{\nu}^{(n-1)}
+30 p_{\nu}  p_{\nu\nu} \tau^2 A_{\nu\nu}^{(n-2)} - 10 p_{\nu\nu} \tau &A_{\nu\nu\nu}^{(n-3)}\Big)
\nonumber\\ 
+ \tau^2\Big(3p_{\nu} \tau A_{\nu}^{(n)} - 3  p_{\nu} \tau A_{\nu\nu}^{(n-1)} + A_{\nu\nu\nu}^{(n-2)} - 3 p_{\nu} p_{\nu\nu} \tau^2 A^{(n)} - 3 p_{\nu\nu} \tau& A_{\nu}^{(n-1)}\Big)\nonumber  \\
+6\tau^4 \sum_{j=0}^{N-1}w^{(j)}A_{\nu}^{(n-2-j)} - c \tau^4 A_{\nu}^{(n-2)}- \nu \tau^4 A_{\nu}^{(n)}& + \tau^5 A_{\tau}^{(n)} = 0.\label{e:2genArecur_kdv5}
\end{align}
Note that there are no terms present due to the $R R_{\nu}$ in this expression, as these only appear at $\mathcal{O}(\e^{-2p\tau/\eps})$. Continuing to expand while taking this term into account would lead to the higher exponentials typically seen in transseries expansions for nonlinear differential equations; we will discuss this in Section \ref{s:end}. 

We may again select the dependence of $A^{(n)}$ on $\tau$ from \eqref{e:f0local}, but some care is required. There are terms in \eqref{e:2genArecur_kdv5} that arise from the nonlinear contribution that do not fit this dependence; however, we will find that these terms do not affect the equation for $f^{(n)}$ at leading order as $n \to \infty$. By substituting \eqref{e:f0local} and retaining only the terms that are dominant in the limit $n \to \infty$ (assuming a factorial-over-power ansatz for $f^{(n)}$ of the form \eqref{e:lateFlocal}), we obtain \eqref{e:frecur}. 

The analysis then proceeds identically to Section \eqref{s:2gen_lkdv5}, obtaining the singulant $\mu$ and prefactor $F$ from \eqref{e:fsingulant} and \eqref{e:fpref} respectively. If we then set \eqref{e:scaling} in the remainder equation \eqref{e:KdV5R_hom}, we obtain an equation similar to that of \eqref{e:2geninnerA}, with additional terms due to the nonlinear contributions. The equation for $A$ is
\begin{align}\nonumber
\Big(3\tau^4 p_{\nu}p_{\nu\nu} + 10 \tau^4 p_{\nu}^3 p_{\nu\nu} A - \tau^4 \nu A_{\nu} + 3 \tau^4 p_{\nu}^2 A_{\nu} + 5 \tau^4 p_{\nu}^4 A_{\nu}\Big)- \eps\Big(15\tau^3 p_{\nu}p_{\nu\nu}A +\tau^3 p_{\nu\nu\nu}A + 10 \tau^3 p_{\nu}^2 p_{\nu\nu\nu}A 
\nonumber&\\
+\,3 \tau^3 p_{\nu\nu}A_{\nu} +30 \tau^3p_{\nu}^2 p_{\nu\nu}A_{\nu} + 3 \tau^3 p_{\nu}A_{\nu\nu} + 10\tau^3 p_{\nu}^3 A_{\nu\nu} \Big)\, +\,\eps^2\Big(5 \tau^2 p_{\nu}p_{\nu\nu\nu\nu}A 
+15 \tau^2 p_{\nu\nu}^2 A_{\nu} 
\nonumber&\\
+ \,20 \tau^2 p_{\nu}p_{\nu\nu\nu}A_{\nu} 
+ 30 \tau^2 p_{\nu}p_{\nu\nu}A_{\nu\nu} + \tau^2 A_{\nu\nu\nu}+ 10 \tau^2 p_{\nu}^2 A_{\nu\nu\nu}\Big) -\eps^3\Big(\tau p_{\nu\nu\nu\nu\nu}A+5 \tau p_{\nu\nu\nu\nu}A_{\nu}
\nonumber&\\
+\, 10 \tau p_{\nu\nu\nu}A_{\nu\nu} + 10 \tau p_{\nu\nu} A_{\nu\nu\nu} - 5 \tau p_{\nu}A_{\nu\nu\nu\nu}\Big) + \eps^4 A_{\nu\nu\nu\nu\nu} - 6\tau^4 A\left( \tau p_{\nu}A + \eps A_{\nu}\right)\frac{\e^{-p\tau/\eps}}{\sqrt{\eps}}
\nonumber&\\
 +\, 6 \eps^2 \tau^4 A \sum_{j=0}^{N-1}\eps^{2j}w_{\nu}^{(j)} 
+ 6\eps \tau^4 \left(\tau  p_{\nu}  A + \eps A_{\nu}\right)\sum_{j=0}^{N-1}\eps^{2j}w^{(j)} 
 - \eps c \tau^4 ( \tau p_{\nu} A + \eps  A_{\nu})+ \tau^5A_{\tau} &= 0.\label{e:2geninnerA_KdV5}
\end{align} 
If we define the inner coordinate $\nu = 2 + \eps^{1/2}\tau^{-1/2}y$ from \eqref{e:2geny} and take the leading-order equation in the limit that $\eps \to 0$, we obtain the same inner expansion as in the analysis of lKdV5 from \eqref{e:2geninner}. From this point, the remaining analysis is identical; there are extra terms on the right-hand side of the equation for the second-generation truncation remainder $Q$, but these terms are asymptotically subdominant in the limit $\eps \to 0$ (and hence $N \to \infty$). Consequently, we find that the asymptotic contribution that appears across that second-generation Stokes line is given by 
\begin{equation}\label{e:newoscillations_KdV5}
Q \sim -\frac{2\pi\i\Lambda}{\eps^{2}} \e^{-\i\nu\tau/\eps} \quad \mathrm{as} \quad \eps \to 0,
\end{equation}
as in the analysis of lKdV5. All that remains is to determine $\Lambda$ by comparing this with the stationary waves \eqref{e:kdv5steadywaves}. In this case, for the solution to be waveless at $\tau = 0$, the precise cancellation requires that $\Lambda \approx 19.97$. 

\subsection{Numerical Comparison}

\begin{figure}
\centering
\includegraphics[width=0.75\textwidth]{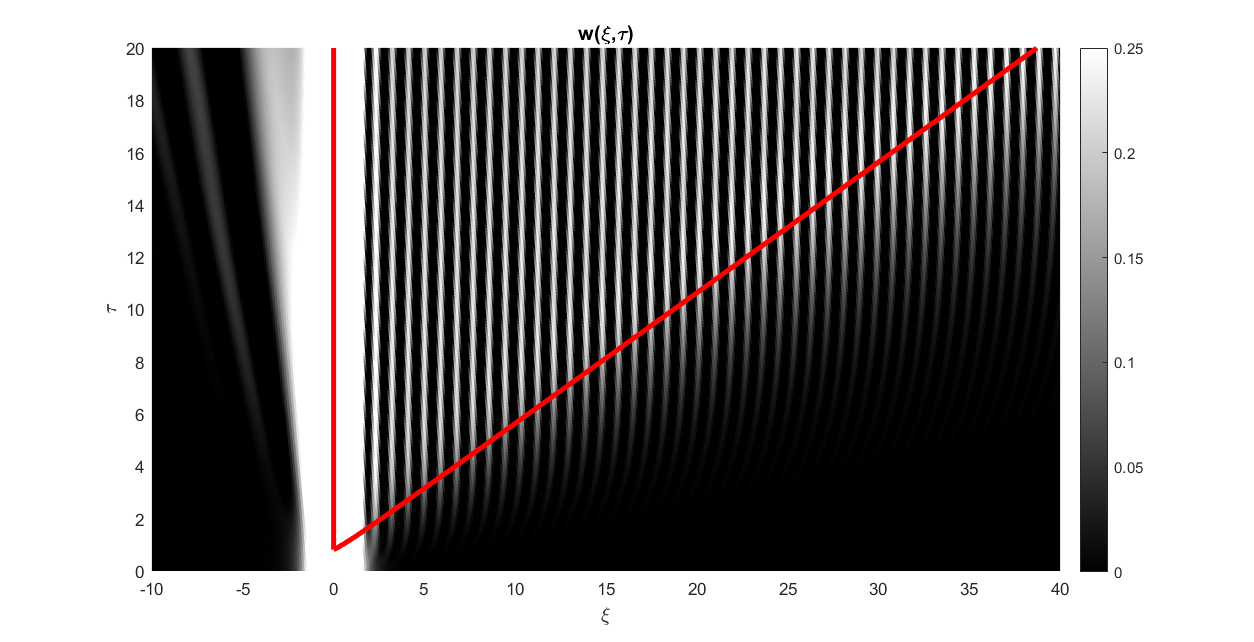}
\caption{Comparison of the numerical solution of KdV5 with $c = 6$ from Figure \ref{fig:KdV5_num} against the Stokes lines obtained by solving \eqref{e:StokesCond} exactly (shown in red). As in the lKdV case, the effect of the preceding transient waves can be seen ahead of the Stokes line, but the waves within the Stokes region are stationary and have constant amplitude.}
\label{fig:KdV5_SL_num}
\end{figure}

In order to validate the asymptotic predictions, we computed the behaviour of the KdV5 system with representative values of $c = 6$ and $\eps = 0.15$ up to $t = 9/20$, which corresponds to $\tau = 20$. We earlier presented this system in Figure \ref{fig:KdV5_num}. We employed a split-step method to compute the solution behaviour (see, for example, \cite{trefethen2000spectral}) with $2^{12}$ points on the periodic domain $\xi \in [-20,120]$. This domain was chosen to be sufficiently large that no visible waves reached the boundaries for $\tau \in [0,20]$. The results of this computation may be seen in Figure \ref{fig:KdV5_SL_num}.

These results indicate that the asymptotic predictions accurately describe the wave propagation within the boundary layer. The transient contribution propogates ahead of the Stokes line, and switches on the stationary ripples that persist as $\tau \to \infty$. These stationary ripples are contained within the Stokes lines, and eventually spread to all $\xi > 0$. The waves that emerge on the left-hand side of the central pulse are algebraically small effects in $\eps$, and are not the consequence of Stokes' phenomenon.

\section{Discussion and Conclusions}\label{s:end}

In this study, we were able to calculate the behaviour of exponentially small oscillations that appear in solutions to a singularly perturbed fifth-order KdV equation (KdV5). We found that the oscillations appear and propagate within a small-time boundary layer in the solution. We constructed an auxiliary linear PDE to determine the appropriate boundary-layer scaling and the algebraic form of the late-order asymptotic terms. This PDE is not a conventional linearisation; instead, it was chosen so that it shares the same late-order behaviour as the full nonlinear KdV5 system. The resulting boundary layer has width $\mathcal{O}(\eps^2)$. This scaling is not predicted by any algebraic balance in $\eps$, but instead emerges only through analysis of the late-order terms and the associated exponentially small contributions.

Within this boundary layer, two exponentially small contributions must be considered: stationary oscillations, whose amplitude remains constant in time, and transient ripples, whose amplitude decays in the far field. Although our primary focus is on the stationary waves, the transient ripples are essential for a complete description of the associated Stokes phenomenon. We found that the stationary oscillations are switched on across a Stokes line generated by the leading-order singularity and switched off across a second-generation Stokes line induced by the transient ripple. Determining the boundaries of the region containing the stationary oscillations shows that they propagate rapidly, with speed $x \sim 2t/\eps^2$ in the asymptotic limit. Furthermore, the transient ripples decay in the outer limit of the boundary layer and therefore cannot be present in the outer solution. This confirms that the non-decaying oscillations persist beyond the boundary layer, which is consistent with the one-sided solution from \cite{Grimshaw,trinh2010exponential, benilov1993generation}. 

Although both factorial-over-power analysis and the method of steepest descents could be applied to the auxiliary linear problem (lKdV5), the latter was not available for the nonlinear KdV5 equation. We therefore applied factorial-over-power analysis directly to the nonlinear system and found that the resulting late-order behaviour and exponentially small contributions closely mirror those of the auxiliary linear problem. Analysis of lKdV5 also identified the Stokes and higher-order Stokes lines relevant to the nonlinear problem, allowing us to focus only on the interactions that contribute to the solution behaviour.

The most challenging aspect of the analysis was determining the inner expansion in the neighbourhood of the singularities. Typically, one computes the inner solution near a singularity, either analytically or numerically, and matches it to the late-order ansatz to determine a multiplicative prefactor. In some linear problems, such as \cite{chapman2005exponential}, this prefactor can be calculated exactly; for nonlinear problems it is usually obtained numerically from the computed series coefficients. For ordinary differential equations, this matching determines a single constant and therefore needs to be performed only once. In the present PDE setting, however, the corresponding prefactor $B$ (or $\tilde{B}$) depends on the characteristic variable $\nu$. A direct numerical approach would therefore require the computation to be repeated across all relevant values of $\nu$. Beyond its computational cost, such an approach obscures the key information: the locations of singularities of $B(\nu)$, their strength, and the local behaviour in their neighbourhoods.

Rather than pursuing this numerical approach, described in Appendix \ref{a:inners}, we exploited information obtained from the analysis of lKdV5. That analysis showed that the stationary waves are confined to a region bounded by a Stokes line describing the interaction between the stationary oscillations and the transient ripple. Identifying the two asymptotic components involved in this interaction was sufficient to determine the value of $\nu$ at which $B(\nu)$ becomes singular. This information is typically available for nonlinear problems, and this method is therefore not limited specifically to the KdV5 equation or even systems with waveless initial conditions. Determining the strength and local behaviour of the singularity required the extra assumption that the oscillations vanished entirely across the second generation Stokes line, and therefore did not appear in the initial condition or persist into the spatial far field. Making this assumption allowed us to identify all of the behaviour required to describe the stationary oscillations, although it was not sufficient to completely determine the amplitude of the transient ripples.

As this is a nonlinear problem, we expect that the full asymptotic solution to the nonlinear KdV5 equation can be expressed as an asymptotic transseries, which is not the case for the linear equation. If a PDE has $N$ singulants, then the general form for the transseries solution is given by
\begin{equation}\label{e:transseries}
w(x,t) \sim \sum_{\mathbf{n} \in \mathbb{N}_0^N} \sigma^\mathbf{n}\e^{-\mathbf{n}\cdot \boldsymbol{\chi}} \Phi_{\mathbf{n}}(x,t,\eps),
\end{equation}
which is indexed by $\mathbf{n} = (n_1, n_2, \ldots, n_N)$ and 
\begin{equation}
\boldsymbol{\chi} = (\chi_1, \chi_2, \ldots, \chi_n), \qquad \sigma^{\mathbf{n}} = \sigma_1^{n_1}\sigma_2^{n_2} \ldots \sigma_N^{n_N}.
\end{equation}
In this expression, the $\sigma_j$ are transseries parameters that change in value across Stokes lines, and $\Phi_{\mathbf n}$ is a Laurent series associated with each combination of exponential contributions. This form contains every possible power and multiplicative combination of the exponential contributions, and therefore admits an infinite hierarchy of Stokes interactions. This structure arises naturally from the remainder equation \eqref{e:KdV5R}. This remainder equation contains a nonlinear term proportional to $R R_{\nu}$, which will force interactions between exponential scales when the transseries ansatz \eqref{e:transseries} is applied. 

If we let $\chi_1 = 0$ (the power series), $\chi_2 = \i\sigma$ (the stationary oscillations) and $\chi_3 = p(\nu)\tau$ (the transient ripples), then the asymptotic contributions studied here correspond to $\mathbf n=(1,0,0,\ldots,0)$, $\mathbf n=(0,1,0,\ldots,0)$, $\mathbf n=(0,0,1,\ldots,0)$ and $\mathbf n=(0,1,1,\ldots,0)$. These elements of the transseries are sufficient to describe the Stokes' phenomenon related to the stationary oscillations at leading order. For more complicated nonlinear problems, higher-order interactions may also play an important role in the relevant solution behaviour.

A further consideration in nonlinear partial differential equations is the role of anti-Stokes lines, on which two distinct exponential contributions are comparable in size. In linear problems, the exponential contributions are independent, and the transseries representation remains valid on anti-Stokes lines. In nonlinear problems, however, contributions of comparable magnitude can interact. Consequently, crossing an anti-Stokes line may alter the dominant balance of the asymptotic solution, and continuation beyond such a line generally requires resummation techniques such as those employed in \cite{lustri2023locating}. These interactions can generate new singular behaviour and increase the complexity of the asymptotic analysis \cite{costin2002formation,mccue2026small}.

This issue does not arise in the present study of KdV5. The asymptotic solution does not contain any anti-Stokes lines in the region where the stationary oscillations and transient ripples are simultaneously active. Consequently, no change in dominant balance occurs within the portion of the solution relevant to the present analysis, and the solution can be continued without the additional resummation procedures required in more general nonlinear problems. Nevertheless, anti-Stokes interactions must typically be taken into account when applying exponential asymptotic methods to nonlinear partial differential equations.

By selecting the initial condition to be the KdV soliton solution, we were able to control the location of singularities in the analytically-continued leading-order solution. For general initial data, the singularity dynamic become far more complicated (see \cite{mccue2026small}), with an infinite number of poles emerging rapidly from points which are singular at $t = 0$. This configuration would produce more complicated dynamics than those seen in the solitary wave problem, as each pole would generate new Stokes lines and interactions.

Finally, an interesting avenue for future work is to investigate whether knowledge of the location and strength of singularities of $B(\nu)$ and $\tilde{B}(\nu)$ can be exploited to reconstruct their local behaviour, and hence find $a^{(0)}$, more efficiently. This could result in a numerical matching procedure that need only be performed once, as in the corresponding theory for ordinary differential equations.

\section{Acknowledgements}

CJL acknowledges the support of Australian Research Council Discovery Project DP240101666.

\section{Conflict of Interest Statement}

The authors have no conflicts of interest to declare. 

\section{References}

\bibliographystyle{plainurl}
\bibliography{reference2}

\appendix
\section{Local Analysis Near Singularities}\label{a:inners}

\subsection{Inner Analysis for lKdV5}

The appropriate neighbourhood is given by $\chi = \mathcal{O}(\eps)$. In \cite{chapman2005exponential}, this motivated the change of variables $\eps z = \chi$, but instead we make use of the result that $\chi = p(\nu)\tau$ and set
\begin{equation}\label{e:innerscale}
y = \nu,\qquad \eps z = \tau, \qquad u = \eps^2 \psi.
\end{equation}
The inner equation becomes
\begin{equation}
\psi_{yyyyy} + z^2 \psi_{yyy} - y z^4 \psi_y + z^5 \psi_z = \eps^2 z^2 \left[\tfrac{c}{2}\sech^2(\tfrac{\sqrt{c}}{2}\eps y z - \tfrac{\i\pi}{\sqrt{2}})\right]_{yyy}.
\end{equation}
Retaining only the leading-order equation as $\eps \to 0$ gives
\begin{equation}
\psi_{yyyyy} + z^2 \psi_{yyy} - y z^4 \psi_y + z^5 \psi_z = \frac{48}{y^5}.
\end{equation}
 The right-hand side simply balances the leading-term in the usual expansion for $\psi$, so we have
 \begin{equation}\label{e:psiinner}
 \psi \sim -\frac{2}{y^2 z^2} + \sum_{n=1}^{\infty} \frac{a^{(n)}(y)}{z^{2n+2}} \quad \as \quad |z| \to \infty.
 \end{equation}
 The recurrence for $a^{(n)}(y)$ is given (for $n \geq 1$, setting $a^{(0)} = -2$ and using the convention that $a^{(k)}(y) = 0$ for $k <0$) by
 \begin{equation}
 a^{(n-2)}_{yyyyy} + a^{(n-1)}_{yyy} - y a^{(n)}_y - 2n a^{(n)} = 0,
 \end{equation}
 with the condition that $\psi$ must decay appropriately as $\tau \to 0$, which will require $a^{(n)} = o(1/y^{2n+2})$. 
Applying this differential-difference equation recursively gives
\begin{align}
a^{(1)}(y) = -\frac{48}{y^5}, \qquad a^{(2)}(y) = -\frac{5040}{y^8} - \frac{1440}{y^7}, \qquad a^{(3)}(y) = -\frac{1209600}{y^{11}} - \frac{725760}{y^{10}},
\end{align}
and this process may be continued indefinitely. To match the inner expansion \eqref{e:psiinner} with the late-order term ansatz \eqref{e:ansatz2}, we determine the inner limit of the late-order terms near the singularity as $\tau \to 0$ and then $\sigma \to 0$. The local late-order terms are
\begin{align}\label{e:ansatz3}
\eps^{2n} u^{(n)}\sim \frac{\eps^{2n} B(\nu)\Gamma(2n + 3/2)}{(4\i \nu/5^{5/4})^{2n+3/2}\tau^{2n+2}}\quad \as \quad n \to \infty, \, \tau \to 0, \, \sigma \to 0.
\end{align}
Matching with the terms of \eqref{e:psiinner}, written in terms of the original coordinates, for large $n$ gives
\begin{equation}\label{e:innereval}
B(\nu) = \lim_{n \to \infty} \frac{a^{(n)}(4\i \nu/5^{5/4})^{2n+3/2}}{\Gamma(2n + 3/2)}.
\end{equation}

\subsection{Inner Analysis for KdV5}\label{s:kdvinner}

The appropriate neighbourhood is identical to that from the lKdV5 analysis, so we set
\begin{equation}\label{e:innerscale2}
y = \nu,\qquad \eps z = \tau, \qquad w = \eps^2 \psi.
\end{equation}
This gives the inner equation
\begin{equation}
\psi_{yyyyy} + z^2 \psi_{yyy} +z^4 \psi \psi_y- y z^4 \psi_y + z^5 \psi_z = 0.
\end{equation}
Applying an expansion as $|z| \to \infty$ in similar fashion to \eqref{e:psiinner} shows that 
 \begin{equation}\label{e:psiinner2}
 \psi \sim -\frac{12}{y^2 z^2} + \sum_{n=1}^{\infty} \frac{a^{(n)}(y)}{z^{2n+2}} \quad \as \quad |z| \to \infty.
 \end{equation}
This is the step in which the nonlinearity in the original PDE alters the analyysis: the recurrence for $a^{(n)}(y)$ for $n \geq 1$ (setting $a^{(0)} = -12$ and using the convention that $a^{(k)}(y) = 0$ for $k <0$) now contains terms arising due to the nonlinear contribution. It has the form
 \begin{equation}
 a^{(n-2)}_{yyyyy} + a^{(n-1)}_{yyy} + \sum_{m=0}^{n-1}a^{(m)}a^{(n-m-1)}_y - y a^{(n)}_y - 2n a^{(n)} = 0,
 \end{equation}
 with the condition that $\psi$ must decay appropriately as $\tau \to 0$, which will require $a^{(n)} = o(1/y^{2n+2})$. This may be used to calculate subsequent values of $a^{(n)}$ and \eqref{e:innereval} can be evaluated as before.

\end{document}